\documentclass{amsproc}
\usepackage{epsfig, a4wide}
\usepackage[latin1]{inputenc}
\usepackage{amsthm, amssymb}
\usepackage{amsfonts, amscd}
\usepackage{multicol}
\usepackage{color}

\newtheorem{theorem}{Theorem}[section]
\newtheorem{lemma}[theorem]{Lemma}

\theoremstyle{definition}

\newtheorem{example}[theorem]{Example}

\theoremstyle{remark}
\newtheorem{remark}[theorem]{Remark}
\newtheorem{proposition}[theorem]{Proposition}
\newtheorem{corollary}[theorem]{Corollary}

\numberwithin{equation}{section}
\usepackage[all]{xy}

\newcommand{\C}{{\mathbb{C}}}
\renewcommand{\H}{{\mathbb{H}}}

\newcommand{\Z}{{\mathbb{Z}}}

\newcommand{\R}{{\mathbb{R}}}
\newcommand{\N}{{\mathbb{N}}}

\begin{document}

\title{Classification of complex projective towers up to dimension $8$ and cohomological rigidity}

\author{Shintar\^o KUROKI}
\address{Osaka City University Advanced Mathematical Institute, Osaka, Japan}
\email{kuroki@scisv.sci.osaka-cu.ac.jp}
\author{DongYoup SUH}
\address{School of Mathematical Science,  Korea Advanced Institute of Science and Technology, Daejeon, Korea}
\email{dysuh@math.kaist.ac.kr}
\thanks{Both authors were supported in part by
the National Research Foundation of Korea(NRF) grant funded by the Korea government(MEST)(No. 2011-0001181).
The first author is partially supported by the JSPS Institutional Program for Young Researcher Overseas Visits
" Promoting international young researchers in mathematics and mathematical sciences led by OCAMI ".}

\subjclass[2010]{}


\begin{abstract}
A complex projective tower or simply a $\mathbb CP$-tower is an iterated complex projective fibrations starting from
a point. In this paper we classify all 6-dimensional $\mathbb CP$-towers up to diffeomorphism, and as a consequence,
we show that all such manifolds are cohomologically rigid, i.e., they are completely determined up to diffeomorphism
by their cohomology rings. 
We also show that cohomological rigidity is not valid for  8-dimensional $\mathbb CP$-towers
by classifying some $\mathbb CP^1$-fibrations over $\mathbb CP^3$ up to diffeomorphism. 
As a corollary we show that such $\mathbb CP$-towers are diffeomorphic if they are homotopy equivalent.
\end{abstract}

\maketitle
\tableofcontents

\section{Introduction}
\label{sect1}

A {\it complex projective tower} (or simply a {\it $\C P$-tower}) of height $m$ is a sequence of complex projective fibrations
\begin{eqnarray}
\label{sequence}
\xymatrix{
& C_{m} \ar[r]^{\pi_{m} \quad } & C_{m-1} \ar[r]^{\pi_{{m-1}}} & \cdots  \ar[r]^{\pi_{2} } & C_{1} \ar[r]^{\pi_{1} \quad \quad \quad } & C_{0}=\{\text{a point}\}
}
\end{eqnarray}
where $C_{i}=P(\xi_{i-1})$ is the projectivization of a complex vector bundle $\xi_{i-1}$ over $C_{i-1}$.
It is also called an {\em $m$-stage $\C P$-tower}.
We call each $C_{i}$ the {\em $i$th stage} of the tower.
Hence a $\C P$-tower is an iterated complex projective bundles starting from a point.

The $\C P$-towers contain many interesting classes of manifolds.
For example, if each complex vector bundle $\xi_{i}$ is a Whitney sum of complex line bundles, such $\C P$-tower is a {\it generalized Bott tower}, introduced in \cite{CMS10}.
If each $\xi_{i}$ is a sum of two complex line bundles, then it is a {\it Bott tower}, introduced in \cite{BoSa} (also see \cite{GrKa}).
In particular, Hirzebruch surfaces are nothing but $2$-stage Bott towers.
Moreover, flag manifolds of type $A$, i.e., $U(n+1)/T^{n+1}\cong \mathcal{F}\ell(\C^{n+1})$, and type $C$, i.e., $Sp(n)/T^{n}$
have $n$-stage $\C P$-tower structures, see Example~\ref{flag manifold} and \ref{flag manifold2}, and the Milnor surface $H_{ij}\subset \C P^{i}\times \C P^{j}$ has a structure of $2$-stage $\C P$-tower, see Example~\ref{Milnor manifold}.

It is well known that there are only two diffeomorphism types of Hirzebruch surfaces, namely, $\C P^{1}\times \C P^{1}$ and $\C P^{2}\#\overline{\C P^{2}}$, and their cohomology rings are not isomorphic.
Hence, Hirzebruch surfaces are classified up to diffeomorphism by their cohomology rings.
One might ask whether the same is true for Bott towers or generalized Bott towers.
Namely, the cohomological rigidity question for (generalized) Bott towers asks whether the diffeomorphism classes of (generalized) Bott towers are determined by their cohomology rings.
There are some partial affirmative answers to the question in \cite{CMS10, CPS, MaPa}, and we refer the reader to \cite{CMS11} for the
summary of the most recent developments about the question.
In particular, the class of $m$-stage Bott towers for $m\le 4$ (\cite{Ch} and \cite{CMS10}) and the class of $2$-stage generalized Bott towers \cite{CMS10} are cohomologically rigid, i.e., their diffeomorphism types are determined by their cohomology rings.

Since the (generalized) Bott tower is a special kind of $\C P$-towers, one might ask the cohomological rigidity question for $\C P$-towers.
On the other hand, if one note that the cohomology ring of a projective bundle $P(\xi)$ is determined by the cohomology ring of the base space of $P(\xi)$ and the Chern classes of the complex bundle $\xi$ (see \eqref{BH-formula}),
then the expectation for the affirmative answer to the question can not be high,
because complex vector bundles are not classified by their Chern classes in general.
Therefore, it might be interesting to determine whether cohomological rigidity indeed fails to hold for $\C P$-towers,
and if so, exactly in what dimension, does it fail?
In this paper, we answer these questions by complete classification of $\C P$-towers up to dimension $6$, and
some special $2$-stage $\C P$-towers of dimension $8$.

We now describe our classification results. Note that the only $2$-dimensional $\C P$-tower is $\C P^{1}$.
Any $4$-dimensional $\C P$-tower is either $\C P^{2}$ or a $2$-stage $\C P$-tower which is in fact nothing but a Hirzebruch surface. So they are
either $H_{0}:=\C P^{1}\times\C P^{1}$ or $H_{1}:=\C P^{2}\#\overline{\C P^{2}}$.
For $6$-dimensional $\C P$-towers, we have to consider one-stage $\C P$-tower which is $\C P^{3}$, two-stage $\C P$-towers, and three-stage $\C P$-towers separately.
For $2$-stage $6$-dimensional $\C P$-towers, there are two cases; the cases when the first stages are $C_{1}=\C P^{1}$ and $C_{1}=\C P^{2}$.
When $C_{1}=\C P^{1}$, then $C_{2}=P(\xi)$ where $\xi$ is a sum of three line bundles.
Therefore, $C_{2}$ must be a $2$-stage generalized Bott tower, which is completely determined in \cite{CMS10}.
In fact, there are only three diffeomorphism types $P(\gamma^{k}_{1}\oplus\epsilon\oplus\epsilon)\to \C P^{1}$ for $k=0,\ 1,\ 2$, where
$\gamma_{1}$ is the tautological line bundle over $\C P^{1}$.

For $2$-stage $6$-dimensional $\C P$-towers with $C_{1}=\C P^{2}$, the second stage $C_{2}=P(\xi)$, where $\xi$ is a rank $2$-complex vector bundle over $\C P^{2}$, which is determined by its Chern classes $c_{1}\in H^{2}(\C P^{2})\simeq \Z$ and $c_{2}\in H^{4}(\C P^{2})\simeq \Z$.
It is proved that the diffeomorphism types of such $\C P$-towers are $P(\eta_{(0,\alpha)})\to \C P^{2}$ and $P(\eta_{(1,\alpha)})\to \C P^{2}$
for $\alpha\in H^{4}(\C P^{2})\simeq \Z$, where $\eta_{(s,\alpha)}$ is a $\C$-vector bundle over $\C P^{2}$ whose Chern classes are $(c_{1},c_{2})=(s,\alpha)$.

For $3$-stage $\C P$-towers $C_{3}\to C_{2}\to C_{1}$, there are two cases, i.e., when $C_{2}=H_{0}=\C P^{1}\times \C P^{1}$ and
$C_{2}=H_{1}=\C P^{2}\#\overline{\C P^{2}}$.
Then $C_{3}=P(\xi)$ where $\xi$ is a complex $2$-dimensional vector bundle over $C_{2}$.
Again, it is proved in Lemma \ref{basic classification} that $\xi$ is classified by its Chern classes $c_{1}$ and $c_{2}$.
Let $\eta_{(s,r,\alpha)}$ (resp.\ $\xi_{(s,r,\alpha)}$) be the complex $2$-dimensional bundle
over $\C P^{1}\times \C P^{1}$ (resp.\ $\C P^{2}\#\overline{\C P^{2}}$) whose first Chern class
$c_{1}(\eta_{(s,r,\alpha)})=(s,r)\in H^{2}(\C P^{1}\times \C P^{1})\simeq \Z\oplus\Z$ (resp.\ $c_{1}(\xi_{(s,r,\alpha)})=(s,r)\in H^{2}(\C P^{2}\#\overline{\C P^{2}})$) and the second Chern class
$c_{2}(\eta_{(s,r,\alpha)})=\alpha\in H^{4}(\C P^{1}\times \C P^{1})\simeq \Z$ (resp.\ $c_{2}(\xi_{(s,r,\alpha)})=\alpha\in H^{4}(\C P^{2}\#\overline{\C P^{2}})$).
Then, it is proved that all diffeomorphism types of  $3$-stage $\C P$-towers are
$P(\zeta_{(s,r,\alpha)})\to H_{0}$ and $P(\xi_{(s,r,\alpha)})\to H_{1}$ for $\alpha\in \Z$ and $(s,r)=(0,0)$, $(1,0)$ or $(1,1)$.

We thus have the following classification result of $6$-dimensional $\C P$-towers.
\begin{theorem}
\label{main}
Any $6$-dimensional $\C P$-tower is diffeomorphic to one of the following distinct manifolds:
\begin{itemize}
\item $\C P^3;$
\item $P(\gamma^{k}_{1}\oplus\epsilon\oplus\epsilon)\to \C P^1$ for $k=0,1,2;$
\item $P(\eta_{(0,\alpha)})\to \C P^2$ for $\alpha\in \Z\setminus \{0\};$
\item $P(\eta_{(1,\alpha)})\to \C P^2$ for $\alpha\in \Z;$
\item $P(\zeta_{(0,0,\alpha)})\to H_{0}$ for $\alpha\in \Z_{\ge 0};$
\item $P(\zeta_{(1,0,\alpha)})\to H_{0}$ for $\alpha\in \Z_{\ge 0};$
\item $P(\zeta_{(1,1,\alpha)})\to H_{0}$ for $\alpha\in \N;$
\item $P(\xi_{(0,0,\alpha)})\to H_{1}$ for $\alpha\in \N;$
\item $P(\xi_{(1,0,\alpha)})\to H_{1}$ for $\alpha\in \Z_{\ge 0};$
\item $P(\xi_{(1,1,\alpha)})\to H_{1}$ for $\alpha\in \Z,$
\end{itemize}
where $H_{0}:=\C P^{1}\times\C P^{1}$, $H_{1}:=\C P^{2}\#\overline{\C P^{2}}$, and the symbols $\N$, $\Z_{\ge 0}$ and $\Z$ represent natural numbers, non-negative integers and integers, respectively.
\end{theorem}
Since the cohomology rings of the manifolds in Theorem \ref{main} are not mutually isomorphic, we have the following corollary
on cohomological rigidity of $\C P$-towers.
\begin{corollary}
\label{main-cor}
Let $M_{1}$ and $M_{2}$ be two $\C P$-towers of dimension less than or equal to $6$.
Then, $M_{1}$ and $M_{2}$ are diffeomorphic if and only if their cohomology rings
$H^{*}(M_{1})$ and $H^{*}(M_{2})$ are isomorphic.
\end{corollary}
This corollary is a generalization of the cohomological rigidity theorem
for Bott manifolds up to dimension less than or equal to $6$ proved in \cite{CMS10}.
Note that cohomology ring does not determine the tower structure of $\C P$-tower.
For example, the trivial $\C P^2$-bundle over $\C P^1$, i.e., $P(\gamma^{0}_{1}\oplus\epsilon\oplus\epsilon)$ and 
the trivial $\C P^1$-bundle over $\C P^2$, i.e., $P(\eta_{(0,0)})$ are diffeomorphic to $\C P^2\times \C P^1$. 
However, they are not equivalent as a fibre bundle.
Namely, their tower structures are different.

For Bott manifolds of dimension $8$ cohomological rigidity theorem is also proved to be true by Choi in \cite{Ch}.
However, it is not the case for $\C P$-towers.
Namely, we classify some special class of $8$-dimensional $2$ stage $\C P$-towers $C_{2}\to C_{1}$ when $C_{1}=\C P^{3}$.
In this case $C_{2}=P(\xi)$ where $\xi$ is a complex $2$-dimensional vector bundle over $\C P^{3}$.
By the result of Atiyah and Rees \cite{AtRe}, any complex $2$-dimensional vector bundle $\xi$ over $\C P^{3}$ is determined by its first and the second Chern classes $c_{1}$ and $c_{2}$ and an
invariant $\alpha\in \Z_{2}$ which is $0$ when $c_{1}$ is odd.
Let $\eta_{(\alpha,c_{1},c_{2})}$ be the complex $2$-dimensional vector bundle with the given invariants $\alpha$, $c_{1}$ and $c_{2}$.
Then we have the following classification theorem of $P(\eta_{(\alpha,c_{1},c_{2})})$.
\begin{theorem}
\label{main-thm2}
Let $M=P(\eta_{(\alpha,c_{1},c_{2})})$ be the projectivization of a $2$-dimensional complex vector bundle $\eta_{(\alpha,c_{1},c_{2})}$ over $\C P^{3}$.
Then the following holds:
\begin{itemize}
\item
if $c_{1}\equiv 1\ ({\rm mod}\ 2)$, then $M$ is diffeomorphic to 
\begin{itemize}
\item $N(t)=P(\eta_{(0,1,t)})$
\end{itemize}
where $t=c_{2}-\frac{c_{1}^{2}-1}{4}\in \Z$;
\item
if $c_{1}\equiv 0\ ({\rm mod}\ 2)$ and $\frac{u(u+1)}{12}\in \Z$ where $u=c_{2}-\frac{c_{1}^{2}}{4}\in \Z$,
then $M$ is diffeomorphic to one of the following distinct manifolds:
\begin{itemize}
\item $M_{0}(u)=P(\eta_{(0,0,u)});$
\item $M_{1}(u)=P(\eta_{(1,0,u)}).$
\end{itemize}
\end{itemize}
\end{theorem}
By the Borel-Hirzebruch formula (\ref{BH-formula}), we have $H^{*}(M_{0}(u))\simeq H^{*}(M_{1}(u))$, while $M_{0}(u)$ 
is not diffeomorphic to $M_{1}(u)$ if $\frac{u(u+1)}{12}\in \Z$.
For example, if $u=0$ then there are just the following two cases: the trivial bundle $M_{0}(0)\simeq \C P^1\times \C P^3$, and the non-trivial bundle $M_{1}(0)$.
This  proves that $8$-dimensional $\C P$-towers are not cohomologically rigid.

On the other hand, we prove that $\pi_{6}(M_{0}(u))\not\simeq \pi_{6}(M_{1}(u))$ if $\frac{u(u+1)}{12}\in \Z$
in Proposition~\ref{homotopy}.
Therefore, we have the following rigidity result.
\begin{corollary}
\label{main-cor2}
Let $\mathcal{C}_{1}$ be the set of $P(\eta_{(\alpha,c_{1},c_{2})})$ with $c_{1}\equiv 1\ ({\rm mod}\ 2)$, 
and $\mathcal{C}_{2}$ be the set of $P(\eta_{(\alpha,c_{1},c_{2})})$ with $c_{1}\equiv 0\ ({\rm mod}\ 2)$ and 
$\frac{u(u+1)}{12}\in \Z$ where $u=c_{2}-\frac{c_{1}^{2}}{4}$.
Then the following holds:
\begin{itemize}
\item $\mathcal{C}_{1}$ satisfies cohomological rigidity, i.e., 
$M_{1}$, $M_{2}\in \mathcal{C}_{1}$ are diffeomorphic if and only if their cohomology rings are isomorphic;
\item $\mathcal{C}_{2}$ satisfies homotopical rigidity, i.e.,
$M_{1}$, $M_{2}\in \mathcal{C}_{2}$ are diffeomorphic if and only if they are homotopic.
\end{itemize}
\end{corollary}


The organization of this paper is as follows.
In Section~\ref{sect2}, we prepare some basics  and  some examples.
In Section~\ref{sect3}, we classify $6$-dimensional $\C P$-towers with height $2$ up to diffeomorphism.
In Section~\ref{sect4}, we classify $6$-dimensional $\C P$-towers with height $3$.
Theorem~\ref{main} is proved as a consequence of the classification. 
In Section~\ref{sect5}, we classify the projectivizations of $2$-dimensional complex vector bundles over $\C P^{3}$, and Theorem~\ref{main-thm2} is proved.

\section{Some preliminaries}
\label{sect2}

In this section, we prepare some basic facts which will be used in later sections.
Let $\xi$ be an $n$-dimensional complex vector bundle over a topological space $X$, and let $P(\xi)$ denote its projectivization.
Then the {\it Borel-Hirzebruch formula} in \cite{BoHi} says
\begin{eqnarray}
\label{BH-formula}
H^{*}(P(\xi);\Z)\simeq H^{*}(X;\Z)[x]/\langle x^{n+1}+\sum_{i=1}^{n}(-1)^{i}c_{i}(\pi^{*}\xi)x^{n+1-i} \rangle
\end{eqnarray}
where
$\pi^{*}\xi$ is the pull-back of $\xi$ along $\pi:P(\xi)\to X$ and
$c_{i}(\pi^{*}\xi)$ is the $i$th Chern class of $\pi^{*}\xi$.
Here $x$ can be viewed as the first Chern class of the canonical line bundle over $P(\xi)$, i.e., 
the complex $1$-dimensional sub-bundle $\gamma_{\xi}$ in $\pi^{*}\xi \to P(\xi)$ such that the restriction $\gamma_{\xi}|_{\pi^{-1}(a)}$ is the canonical line bundle over $\pi^{-1}(a)\cong \C P^{n-1}$ for all $a\in X$.
Therefore $\deg x=2$.
Since it is well-known that the induced homomorphism $\pi^{*}:H^{*}(X;\Z)\to H^{*}(P(\xi);\Z)$ is injective,
we often confuse $c_{i}(\pi^{*}\xi)$  with $c_{i}(\xi)$.

We apply the formula~\eqref{BH-formula} to an $m$-stage $\C P$-tower
\begin{eqnarray*}
\xymatrix{
& C_{m} \ar[r]^{\pi_{m} \quad } & C_{m-1} \ar[r]^{\pi_{{m-1}}} & \cdots  \ar[r]^{\pi_{2} } & C_{1} \ar[r]^{\pi_{1} \quad \quad \quad } & C_{0}=\{\text{a point}\}
}
\end{eqnarray*}
with  $C_i=P(\xi_{i-1})$, to get the following isomorphisms.
\begin{eqnarray}
H^{*}(C_{m};\Z)&\simeq & H^{*}(C_{m-1};\Z)[x_{m}]/\langle x_{m}^{n_{m}+1}+\sum_{i=1}^{n_{m}}(-1)^{i}c_{i}(\xi_{m-1})x_{m}^{n_{m}+1-i} \rangle \nonumber  \\
\nonumber &\simeq & H^{*}(C_{m-2};\Z)[x_{m-1}, x_{m}]/\langle
x_{k}^{n_{k}+1}+\sum_{i=1}^{n_{k}}(-1)^{i}c_{i}(\xi_{k})x_{k}^{n_{k}+1-i}\ |\ k=m-1,\ m \rangle \nonumber  \\
 & \vdots & \nonumber \\
 &\simeq &\Z[x_{1},\ldots, x_{m}]/\langle
x_{k}^{n_{k}+1}+\sum_{i=1}^{n_{k}}(-1)^{i}c_{i}(\xi_{k})x_{k}^{n_{k}+1-i}\ |\ k=1,\cdots, m \rangle. \label{cohom-ring}
\end{eqnarray}

In order to prove the main theorem, we often use the following lemmas.
\begin{lemma}
\label{tensor line budle}
Let $\gamma$ be any line bundle over $M$,
and let $P(\xi)$ be the projectivization of a complex vector bundle $\xi$ over $M$.
Then, $P(\xi)$ is diffeomorphic to $P(\xi\otimes \gamma)$.
\end{lemma}
\begin{proof}
By the definition of the projectivization of a complex vector bundle, the statement follows immediately.
\end{proof}

\begin{lemma}
\label{formula of Chern class}
Let $\gamma$ be a complex line bundle, and let $\xi$ be a $2$-dimensional complex vector bundle over a manifold $M$.
Then the Chern classes of the  tensor product $\xi\otimes\gamma$ are as follows.
\begin{eqnarray*}
& &c_{1}(\xi\otimes\gamma)=c_{1}(\xi)+2c_{1}(\gamma); \\
& &c_{2}(\xi\otimes\gamma)=c_{1}(\gamma)^2+c_{1}(\gamma)c_{1}(\xi)+c_{2}(\xi).
\end{eqnarray*}
\end{lemma}
\begin{proof}
Let us consider the following pull-back diagram:
\begin{eqnarray*}
\xymatrix{
& \pi^{*}\xi\otimes\pi^{*}\gamma \ar[r] \ar[d] & \xi\otimes\gamma \ar[d] \\
& P(\xi\otimes\gamma) \ar[r]^{\pi} & M
}
\end{eqnarray*}
Let $\varphi : P(\xi\otimes \gamma)\to P(\xi)$ be the diffeomorphism from Lemma~\ref{tensor line budle}, and let
$\pi_{\xi}: P(\xi)\to M$ be the projection of the fibration. Then we can see easily that $\pi=\pi_{\xi}\circ \varphi$.
Taking the canonical line bundle $\gamma_{\xi}$ in $\pi_{\xi}^{*}\xi$,
we may regard $\pi_{\xi}^{*}\xi\equiv \gamma_{\xi}\oplus\gamma_{\xi}^{\perp}$, where $\gamma_{\xi}^{\perp}$ is the normal (line) bundle of $\gamma_{\xi}$ in $\pi_{\xi}^{*}\xi$.
By using the decomposition $\pi=\pi_{\xi}\circ \varphi$,
we have the following equation:
\begin{eqnarray*}
\pi^{*}c(\xi\otimes\gamma)&=&c(\varphi^{*}\gamma_{\xi}\otimes\pi^{*}(\gamma))c(\varphi^{*}\gamma_{\xi}^{\perp}\otimes\pi^{*}(\gamma)) \\
&=&(1+\varphi^{*}c_{1}(\gamma_{\xi})+\pi^{*}c_{1}(\gamma))(1+\varphi^{*}c_{1}(\gamma_{\xi}^{\perp})+\pi^{*}c_{1}(\gamma)).
\end{eqnarray*}
Because $\pi^{*}c_{1}(\xi)=\varphi^{*}c_{1}(\gamma_{\xi})+\varphi^{*}c_{1}(\gamma_{\xi}^{\perp})$ and $\pi^{*}c_{2}(\xi)=\varphi^{*}c_{1}(\gamma_{\xi})\varphi^{*}c_{1}(\gamma_{\xi}^{\perp})$,
we have
\begin{eqnarray*}
\pi^{*}c_{1}(\xi\otimes\gamma)&=& \pi^{*}c_{1}(\xi)+2\pi^{*}c_{1}(\gamma); \\
\pi^{*}c_{2}(\xi\otimes\gamma)&=& \pi^{*}c_{2}(\xi)+\pi^{*}c_{1}(\xi)\pi^{*}c_{1}(\gamma)+\pi^{*}c_{1}(\gamma)^{2}.
\end{eqnarray*}
As is well-known, $\pi^{*}:H^{*}(M)\to H^{*}(P(\xi\otimes\gamma))$ is injective. Hence
we have the formula in the lemma.
\end{proof}

We now give two examples of $\C P$-towers.
\begin{example}
\label{flag manifold}
The flag manifold $\mathcal{F}l(\C^{n+1})=\{\{0\}\subset V_{1}\subset \cdots \subset V_{n}\subset \C^{n+1}\}$, called type $A$,
is well-known to be diffeomorphic to the homogeneous space $U(n+1)/T^{n+1}(\cong SU(n+1)/T^{n})$.
We will show that the flag manifold $U(n+1)/T^{n+1}$ is a  
$\C P$-tower with height $n$.
Recall that if  $M$ is a smooth manifold with free $K$ action and $H$ is a subgroup of $K$,
then we have a diffeomorphism
$M/H\cong M\times_{K}(K/H) $.
Also recall  that $\C P^{n}\cong U(n+1)/(T^{1}\times U(n))$.
By using these facts, it is easy to check that there is the following $\C P$-tower structure of height $n$  in $U(n+1)/T^{n+1}$:
\begin{eqnarray*}
\begin{array}{c}
U(n+1)\times_{(T^{1}\times U(n))}(U(n)\times_{(T^{1}\times U(n-1))}(U(n-1)\times_{(T^{1}\times U(n-2))}\cdots (U(3)\times_{(T^{1}\times U(2))}\C P^1)\cdots ) \\
\downarrow \\
\vdots \\
\downarrow \\
U(n+1)\times_{(T^{1}\times U(n))}(U(n)\times_{(T^{1}\times U(n-1))}\C P^{n-2}) \\
\downarrow  \\
U(n+1)\times_{(T^{1}\times U(n))}\C P^{n-1} \\
\downarrow  \\
\C P^{n},
\end{array}
\end{eqnarray*}
where the $U(k)$ action on $\C P^{k-1}$ in each stage is induced from the usual $U(k)$ action on $\C ^k$.
\end{example}

\begin{example}
\label{flag manifold2}
The flag manifold of type C is defined by the homogeneous space $Sp(n)/T^{n}$.
We claim that $Sp(n)/T^{n}$ is a $\C P$-tower with height $n$.
It is well known that $Sp(n)/(T^{1}\times Sp(n-1))\cong S^{4n-1}/T^{1}\cong \C P^{2n-1}$, because $Sp(n)/Sp(n-1)\cong S^{4n-1}$.
By using this fact and the method similar to that demonstrated in Example~\ref{flag manifold}, 
it is easy to check that there is the following $\C P$-tower structure of height $n$  in $Sp(n)/T^{n}$:
\begin{eqnarray*}
\begin{array}{c}
Sp(n)\times_{(T^{1}\times Sp(n-1))}(Sp(n-1)\times_{(T^{1}\times Sp(n-2))}\cdots (Sp(2)\times_{(T^{1}\times Sp(1))}\C P^1)\cdots ) \\
\downarrow \\
\vdots \\
\downarrow \\
Sp(n)\times_{(T^{1}\times Sp(n-1))}(Sp(n-1)\times_{(T^{1}\times Sp(n-2))}\C P^{2n-5}) \\
\downarrow  \\
Sp(n)\times_{(T^{1}\times Sp(n-1))}\C P^{2n-3} \\
\downarrow  \\
\C P^{2n-1},
\end{array}
\end{eqnarray*}
where the $Sp(k)$-action on $\C P^{2k-1}$ in each stage is induced from the $Sp(k)$-action on $\C^{2k}(\simeq \H^{k})$ induced by the following representation to $U(2k)$:
\begin{eqnarray*}
A+Bj\longrightarrow 
\left(
\begin{array}{cc}
A & -B \\
\overline{B} & \overline{A}
\end{array}
\right).
\end{eqnarray*}
Here $A,\ B\in M(k;\C)$ satisfy $A\overline{A}+B\overline{B}=I_{k}$ and $BA-AB=O$.
\end{example}
\begin{remark}
By computing the generators of flag manifolds of other types ($B_{n}$ ($n\ge 3$), $D_{n}$ ($n\ge 4$), $G_{2}$, $F_{4}$, $E_{6},\ E_{7},\ E_{8}$),
they do not admit the structure of $\C P$-towers, see \cite{Bo} (or \cite{FIM} for classical types).
\end{remark}

\begin{example}
\label{Milnor manifold}
The {\it Milnor hypersurface} $H_{i,j}\subset \C P^{i}\times \C P^{j}$, $1\le i\le j$ is defined by the following equation (see \cite[Example 5.39]{BuPa}):
\begin{eqnarray*}
H_{i,j}=\{[z_{0}:\cdots : z_{i}]\times [w_{0}:\cdots : w_{j}]\in \C P^{i}\times \C P^{j}\ |\ \sum_{q=0}^{i}z_{q}w_{q}=0\}.
\end{eqnarray*}
We can  show easily that the natural projection onto the first coordinate of $H_{i,j}$ gives the structure of a $\C P^{j-1}$-bundle over $\C P^{i}$.
Moreover, by the proof in \cite[Theorem 5.39]{BuPa}, this bundle may be regarded as the projectivization of $\gamma^{\perp}\subset \epsilon^{j+1}$, where $\epsilon^{j+1}$ is the trivial $\C^{j+1}$-bundle over $\C P^{i}$ and $\gamma^{\perp}$ is the normal bundle of
the canonical line bundle $\gamma$ over $\C P^{i}$ in $\epsilon^{j+1}$.
Therefore, the Milnor hypersurface admits the structure of a $\C P$-tower with height $2$.
\end{example}

\begin{remark}
\label{two-examples}
As is well-known,
both of the flag manifold $U(n+1)/T^{n+1}$ (and $Sp(n)/T^{n}$) with $n\ge 2$ and the Milnor hypersurface $H_{i,j}$ with $i\ge 2$ do not admit the structure of a {\it toric manifold} (see e.g. \cite{BuPa}).
On the other hand, $U(2)/T^{2}\cong Sp(1)/T^{1}\cong \C P^1$ and $H_{1,j}\to \C P^{1}$ are toric manifolds.
\end{remark}

\section{$6$-dimensional $\C P$-towers of height $2$}
\label{sect3}

Let $M$ be a $6$-dimensional $\C P$-tower.
Then, the height of $M$ is at most $3$.
If its height is one, then $M$ is diffeomorphic to $\C P^3$.
Therefore, it is enough to analyze the case when the height is $2$ and $3$.
In this section, we focus on the classification of $6$-dimensional $\C P$-towers of height $2$.

To state the main theorem of this section, we first set up some notation.
Let  $\mathcal{M}_{2}^{6}$ be the set of all $6$-dimensional $\C P$-towers of height $2$, up to diffeomorphisms.
Let $\gamma_{i}$ denote the tautological line bundle over $\C P^i$,
and let $x$ denote the generator $-c_{1}(\gamma_{2})\in H^{2}(\C P^2)$.
Let $\eta_{(s,\alpha)}$ as the complex $2$-dimensional vector bundle over $\C P^{2}$ whose
total Chern class is $1+sx+\alpha x^2$
for $s,\ \alpha\in \Z$, let $P(\eta_{(s,\alpha)})$ be its projectivization.
We now state the main theorem of this section.
\begin{theorem}
\label{2-stage}
The set $\mathcal{M}_{2}^{6}$ consists of the following distinct $\C P$-towers.
\begin{eqnarray*}
& &P(\gamma_{1}\oplus\epsilon\oplus\epsilon)\longrightarrow \C P^1; \\
& &P(\gamma_{1}^{2}\oplus\epsilon\oplus\epsilon)\longrightarrow \C P^1, \ {\textrm where }\    \gamma_{1}^{2}\equiv \gamma_{1}\otimes \gamma_{1}; \\
& &P(\eta_{(0,\alpha)})\longrightarrow \C P^2 \quad for\quad \alpha\in \Z; \\
& &P(\eta_{(1,\beta)})\longrightarrow \C P^2 \quad for\quad \beta\in \Z.
\end{eqnarray*}
\end{theorem}

\begin{proof}
Take $M\in \mathcal{M}_{2}^{6}$.
Then the first stage $C_{1}$  of $M$ is either $\C P^1$ or $\C P^2$.
We treat these two cases separately below.

\subsection*{CASE I: $C_{1}=\C P^1$}

Note that any complex vector bundles over $\C P^1$ decomposes into a Whitney sum of line bundles.
Therefore a $\C P$-tower $M\in \mathcal{M}_{2}^{6}$ with $C_{1}=\C P^1$ is a $2$-stage generalized Bott tower, and such
Bott towers are completely classified in \cite{CMS10}. (See also \cite{CPS}.)
Due to the cited result, we have the following proposition.
\begin{proposition}
\label{rank3}
Let $M\in \mathcal{M}_{2}^{6}$ be a generalized Bot manifold with $C_{1}=\C P^1$.
Then $M$ is diffeomorphic to one of the following three  distinct manifolds:
\begin{eqnarray*}
& & P(\gamma_{1}^{0}\oplus\epsilon\oplus\epsilon)\cong \C P^1\times \C P^2, \ {\textrm where}\ \gamma_{1}^{0}\equiv \epsilon; \\
& & P(\gamma_{1}\oplus\epsilon\oplus\epsilon); \\
& & P(\gamma_{1}^{2}\oplus\epsilon\oplus\epsilon).
\end{eqnarray*}
\end{proposition}

\subsection*{CASE II: $C_{1}=\C P^2$}

Because $\dim M=6$ and $C_{1}=\C P^2$, the bundle $E_{1}\to C_{1}$ is a complex $2$-dimensional vector bundle.
Such vector bundles are determined by their Chern classes $c_{1}$ and $c_{2}$ (see \cite{Sh, Sw}).
Hence, by Lemmas \ref{tensor line budle} and \ref{formula of Chern class},
we may denote $E_{1}$ by $\eta_{(s,\alpha)}$ such that
$c_{1}(\eta_{(s,\alpha)})=s x$  for $s=0,1$ and $c_{2}(\eta_{(s,\alpha)})=\alpha x^2 \in H^{4}(\C P^2)$  for $\alpha\in \Z$.
In Case II, we have the following classification result.
\begin{proposition}
\label{rank2}
The following  are equivalent for $s_{1},\ s_{2}\in \{0,\ 1\}$ and $\alpha_{1},\ \alpha_{2}\in \Z$.
\begin{enumerate}
\item $(s_{1}, \alpha_{1})=(s_{2}, \alpha_{2})$.
\item Two manifolds $P(\eta_{(s_{1}, \alpha_{1})})$ and  $P(\eta_{(s_{2}, \alpha_{2})})$ are diffeomorphic.
\item Two cohomology rings $H^{*}(P(\eta_{(s_{1}, \alpha_{1})}))$ and $H^{*}(P(\eta_{(s_{2}, \alpha_{2})}))$ isomorphic.
\end{enumerate}
\end{proposition}
Theorem \ref{2-stage} follows from Proposition \ref{rank3} and \ref{rank2}.
\end{proof}
It remains to prove Proposition \ref{rank2}.

\begin{proof}[Proof of Proposition \ref{rank2}]
$(1)\Rightarrow (2)$ and $(2)\Rightarrow (3)$ are obvious.
We now prove $(3)\Rightarrow (1)$. We prove this  by proving the three claims:
(1) $H^{*}(P(\eta_{(0, \alpha)}))\not\simeq H^{*}(P(\eta_{(1, \beta)}))$ for every $\alpha,\ \beta\in \Z$,
(2) if $H^{*}(P(\eta_{(0, \alpha_{1})}))\simeq H^{*}(P(\eta_{(0, \alpha_{2})}))$ then $\alpha_{1}=\alpha_{2}$, and
(3) if $H^{*}(P(\eta_{(1, \beta_{1})}))\simeq H^{*}(P(\eta_{(1, \beta_{2})}))$ then $\beta_{1}=\beta_{2}$.

\textbf{Claim 1: $H^{*}(P(\eta_{(0, \alpha)}))\not\simeq H^{*}(P(\eta_{(1, \beta)}))$ for every $\alpha,\ \beta\in \Z$.}
By using the Borel-Hirzebruch formula (\ref{BH-formula}), we have the following isomorphisms:
\begin{eqnarray*}
H^{*}(P(\eta_{(0,\alpha)})) &\simeq & \Z[X,Y]/\langle X^{3}, Y^{2}+\alpha X^{2} \rangle; \\
H^{*}(P(\eta_{(1,\beta)})) &\simeq & \Z[x,y]/\langle x^{3}, y^{2}+xy+\beta x^{2} \rangle,
\end{eqnarray*}
where $\deg X=\deg Y=\deg x=\deg y=2$.
We write the $\Z$-module structures of $H^{*}(P(\eta_{(0,\alpha)}))$ and $H^{*}(P(\eta_{(1,\beta)}))$ by indicating  their generators as follows:
\begin{eqnarray*}
& &\Z\oplus \Z X\oplus\Z Y \oplus \Z X^2\oplus \Z XY \oplus \Z X^2 Y; \\
& &\Z\oplus \Z x\oplus\Z y \oplus \Z x^2\oplus \Z xy \oplus \Z x^2 y.
\end{eqnarray*}

If there exits a graded ring isomorphism $f:H^{*}(P(\eta_{(0,\alpha)}))\to H^{*}(P(\eta_{(1,\beta)}))$, then
we may put $f(X)=ax+by$ and $f(Y)=cx+dy$ for some $a,b,c,d\in \Z$ such that
\begin{eqnarray}
\label{det}
ad-bc=\pm 1.
\end{eqnarray}
Because $f$ preserves the ring structure, we have
\begin{eqnarray*}
f(X^3)&=&(ax+by)^3 \\
&= & (3a^2b -3ab^2+b^3-\beta b^3)x^2y=0; \\
f(Y^2+\alpha X^2)&=&(cx+dy)^2+\alpha (ax+by)^2 \\
&= & (c^2+\alpha a^2-\beta d^2-\alpha \beta b^2)x^2+(2cd+2\alpha ab-d^2-\alpha b^2)xy=0.
\end{eqnarray*}
This implies the following equations:
\begin{eqnarray}
& & \label{from x^3} b(3a^2-3ab+b^2-\beta b^2)=0; \\
& & \label{from y^2-1}  c^2+\alpha a^2-\beta d^2-\alpha \beta b^2=0; \\
& & \label{from y^2-2}  2cd+2\alpha ab-d^2-\alpha b^2=0.
\end{eqnarray}
If $b=0$, then $2c=d=\pm 1$ by \eqref{det} and \eqref{from y^2-2}.
But this contradicts to the fact that $c$ is an integer (i.e., $c\in \Z$).
Hence $b\not=0$,  and by  \eqref{from x^3} we have $3a^2-3ab+b^2-\beta b^2=0$.
We also have the following commutative diagram of free $\Z$-modules.
\begin{eqnarray*}
\xymatrix{
\Z X\oplus \Z Y \ar[d]^{f} \ar[r]^{\cdot X} & \Z X^2\oplus \Z XY \ar[d]^{f} \\
\Z x\oplus \Z y \ar[r]^{\cdot f(X)} & \Z x^2\oplus \Z xy
},
\end{eqnarray*}
where the horizontal maps are induced from the multiplication by $X$ and $f(X)$, respectively.
Let us represent the linear map $\cdot f(X)=\cdot (ax+by): \Z x\oplus \Z y\to \Z x^2\oplus \Z xy$ by the  matrix
\begin{eqnarray*}
A=\left(
\begin{array}{cc}
a & -\beta b \\
b & a-b
\end{array}
\right)
\end{eqnarray*}
with respect to the generators.
Note that  $\cdot X:\Z X\oplus \Z Y\to \Z X^2\oplus \Z XY$ is an isomorphism.
Therefore $\cdot f(X)$ is also an isomorphism, and hence
\begin{eqnarray}
\label{from matrix}
\det A=a^2-ab+\beta b^2=\pm 1.
\end{eqnarray}

Because $b\not=0$, it follows from \eqref{from x^3} and \eqref{from matrix} that we have $b=\pm 1$, $\beta=1$ and $a=0$ or $b$.
If $a=b$, then $c=d$ or $c=-d$ by \eqref{from y^2-1}.
However, it is easy to check that both of these cases give contradictions to \eqref{det} and $c,d\in \Z$.
Hence, $a=0$.
In this case, $\alpha=c^2-d^2$ by \eqref{from y^2-1} and $\alpha=2cd-d^2$ by \eqref{from y^2-2}.
Therefore we have $c=0$ or $2d$.
However, both of these cases give contradictions to \eqref{det} and $c,d\in \Z$.
This establishes that there is no ring isomorphism between $H^{*}(P(\eta_{(0,\alpha)}))$ and $H^{*}(P(\eta_{(1,\beta)}))$.

\textbf{Claim 2: If $H^{*}(P(\eta_{(0, \alpha_{1})}))\simeq H^{*}(P(\eta_{(0, \alpha_{2})}))$, then $\alpha_{1}=\alpha_{2}$.}
By (\ref{BH-formula}), we have the isomorphisms
\begin{eqnarray*}
H^{*}(P(\eta_{(0,\alpha_1)})) &\simeq & \Z[X,Y]/\langle X^{3}, Y^{2}+\alpha_1 X^{2} \rangle, \ {\rm  and}\\
H^{*}(P(\eta_{(0,\alpha_2)})) &\simeq & \Z[x,y]/\langle x^{3}, y^{2}+\alpha_2 x^{2} \rangle.
\end{eqnarray*}
Assume that there exists an isomorphism $f:H^{*}(P(\eta_{(0,\alpha_{1})}))\to H^{*}(P(\eta_{(0,\alpha_{2})}))$ for some $\alpha_{1},\ \alpha_{2}\in \Z$, and let $f(X)=ax+by$ and $f(Y)=cx+dy$, so that $ad-bc=\pm 1$.
Because $f(X^3)=(ax+by)^3=0$, we have that
\begin{eqnarray*}
b(3a^2-b^2\alpha_{2})=0.
\end{eqnarray*}
Suppose $b\not=0$. Then $3a^2-b^2\alpha_{2}=0$.
Because the  map
\begin{eqnarray*}
f:H^{6}(P(\eta_{(0,\alpha_{1})}))=\Z X^2Y\longrightarrow \Z x^2y=H^{6}(P(\eta_{(0,\alpha_{2})})),
\end{eqnarray*}
is an isomorphism, we have
\begin{eqnarray}
\label{top isom}
f(X^2Y)=(ax+by)^2(cx+dy)=\pm x^2y.
\end{eqnarray}
Using \eqref{top isom} and the ring structures, we have that
\begin{eqnarray*}
a^2d+2abc-b^2d\alpha_{2}=\pm 1.
\end{eqnarray*}
Because $3a^2-b^2\alpha_{2}=0$, we have $-2a^2d+2abc=-2a(ad-bc)=\pm 1$.
However, this gives a contradiction to $a\in \Z$, because $ad-bc=\pm 1$.
Hence, $b=0$ and $ad=\pm 1$; in particular, we have $a,\ d=\pm 1$.
Then, we have the following equations:
\begin{eqnarray*}
f(Y^2+\alpha_{1} X^2)&=&(cx+dy)^2+\alpha_{1}(ax+by)^2 \\
&=&(c^2-\alpha_{2}+\alpha_{1}) x^2+2cd xy=0
\end{eqnarray*}
Therefore, we have that $c=0$ and $\alpha_{1}=\alpha_{2}$.
This proves the claim.

\textbf{Claim 3: If $H^{*}(P(\eta_{(1, \beta_{1})}))\simeq H^{*}(P(\eta_{(1, \beta_{2})}))$, then $\beta_{1}=\beta_{2}$.}
By (\ref{BH-formula}), we have the  isomorphisms
\begin{eqnarray*}
H^{*}(P(\eta_{(1,\beta_1)})) &\simeq & \Z[X,Y]/\langle X^{3}, Y^{2}+XY+\beta_1 X^{2} \rangle, \ {\rm and} \\
H^{*}(P(\eta_{(1,\beta_2)})) &\simeq & \Z[x,y]/\langle x^{3}, y^{2}+xy+\beta_2 x^{2} \rangle.
\end{eqnarray*}
Assume that there exists an isomorphism $f:H^{*}(P(\eta_{(1,\beta_{1})}))\to H^{*}(P(\eta_{(1,\beta_{2})}))$ for some $\beta_{1},\ \beta_{2}\in \Z$, and  let $f(X)=ax+by$ and $f(Y)=cx+dy$, so that $ad-bc=\pm 1$.
Because of the relations $f(X^3)=(ax+by)^3=0$ and $f(Y^2+XY+\beta_{1}X^2)=(cx+dy)^2+(ax+by)(cx+dy)+\beta_{1}(ax+by)^2=0$,
we have that
\begin{eqnarray}
\label{from rel of x^3}
&b(3a^2-3ab+b^2-b^2\beta_{2})=0; \\
\label{from rel other1}
&c^2-d^2\beta_{2}+ac-bd\beta_{2}+a^2\beta_{1}-b^2\beta_{1}\beta_{2}=0; \\
\label{from rel other2}
&2cd-d^2+ad+bc-bd+2\beta_{1}ab-\beta_{1}b^2=0.
\end{eqnarray}

We first assume $b=0$.
From the equation $ad-bc=\pm 1$, we have  $a,\ d=\pm 1$.
Now plug $b=0$ and $d=\pm 1$ into \eqref{from rel other2} to get the equation
\begin{eqnarray*}
2c+a=d=\pm 1.
\end{eqnarray*}
Together with $a=\pm 1$, this equation implies that either $c=0$ and  $a=d$, or $c\not=0$ and $c=-a=d$.
Now plug these into \eqref{from rel other1} to obtain $\beta_{1}=\beta_{2}$ in either cases, which proves the claim
when $b=0$.

We now assume $b\not=0$. Then from \eqref{from rel of x^3}, we have $3a^2-3ab+b^2-b^2\beta_{2}=0$.
By using the same  argument as the one used to get \eqref{from matrix}, we have
\begin{eqnarray}
\label{for a}
a^2-ab+\beta_{2} b^2=\epsilon,
\end{eqnarray}
where $\epsilon=\pm 1$.
Substitute \eqref{for a} into the equation $3a^2-3ab+b^2-b^2\beta_{2}=0$.
Then, we obtain the equation
\begin{eqnarray*}
b^2(4\beta_{2}-1)=3\epsilon.
\end{eqnarray*}
Therefore, $b=\pm 1$ and $\beta_{2}=\epsilon=1$.
Hence, together with \eqref{for a}, we have that $a=0$ or $a=b$.

If $a=0$, then $c=\pm 1$ by the equation $ad-bc=\pm 1$.
Substitute these equations into \eqref{from rel other1} and \eqref{from rel other2}.
Then, we have the equations
\begin{eqnarray*}
\beta_{1}=1-d^2-bd=2cd-d^2+bc-bd.
\end{eqnarray*}
Therefore, we have that $(2d+b)c=1$.
Moreover,
because $c=\pm 1$ and $b=\pm 1$, we have $(b,d)=(c,0)$ or $(-c,c)$.
Hence, $\beta_{1}=1=\beta_{2}$.

If $a=b=\pm 1$, then $d-c=\pm 1$ by the equation $ad-bc=\pm 1$.
Put $a=b=\pm 1$ in \eqref{from rel other2} to obtain the equation
\begin{eqnarray}
\label{for alpha}
\beta_{1}=d^2-2cd-bc.
\end{eqnarray}
Moreover, by substituting $a=b=\pm 1$ and $\beta_{2}=1$ into \eqref{from rel other1}, we have
\begin{eqnarray*}
(c-d)(a+c+d)=0.
\end{eqnarray*}
This together with $d-c=\pm 1$ implies that $c+d=-a=\pm 1$.
It follows that either $d=0$ and $c=-a=-b$, or $d=-a=-b$ and $c=0$.
By \eqref{for alpha}, we have $\beta_{1}=1=\beta_{2}$.  This proves the claim, and hence
the proof of the proposition is complete.
\end{proof}

We can show easily  that $P(\eta_{(s,\alpha)})$ is diffeomorphic to $\C P^1 \times \C P^2$ if and only if
$(s,\alpha)=(0,0)$ by comparing their cohomology rings.
Therefore, by Propositions \ref{rank3} and \ref{rank2}, we have Theorem \ref{2-stage}.
Moreover, by Theorem \ref{2-stage}, we have the following corollary.
\begin{corollary}
\label{cohom. rigid. up to 2-stage}
Let $\mathcal{M}_{\le 2}^{6}$ be the class of all $6$-dimensional $\C P$-towers of height at most $2$, up to diffeomorphism.
Then two $\C P$-towers $M$ and $ M'$ in  $\mathcal{M}_{\le 2}^{6}$ are diffeomorphic if and  only if their cohomology rings
$H^{*}(M)$ and $H^{*}(M')$ are isomorphic.
In other words, the class
$\mathcal{M}_{\le 2}^{6}$ is  cohomologically rigid.
\end{corollary}

\section{$3$-stage $6$-dimensional $\C P$-towers}
\label{sect4}

In this section, we focus on  $6$-dimensional $\C P$-towers of height $3$.
The $3$-stage $6$-dimensional $\C P$-towers are of the form
\begin{eqnarray*}
\xymatrix{
P(\xi)\ar[r]^{\C P^1} & H_{k}\ar[r]^{\C P^1} & \C P^1.}
\end{eqnarray*}
Here, $\xi$ is a complex $2$-dimensional vector bundle over $H_{k}$, and  $H_{k}$ is the {\it Hirzebruch surface} $P(\gamma^{k}_{1}\oplus\epsilon)$
where $\epsilon$ is the trivial complex line bundle and
$\gamma_{1}^{k}$ is the $k$-th tensor power of the tautological line bundle $\gamma_{1}$ over $\C P^1$.
As is well known,  $H_{k}$ is diffeomorphic to $H_0$ if $k$ is even, and to $H_1$ if $k$ is odd (see \cite{Hi, MaSu}).

\begin{lemma}
\label{basic classification}
Let ${\rm Vect}_{\C}^{2}(H_{k})$ be the set of complex $2$-dimensional  vector bundles over $H_{k}$ up to isomorphisms.
Then the correspondence
\begin{eqnarray*}
\begin{array}{ccc}
{\rm Vect}_{\C}^{2}(H_{k}) & \stackrel{c}{\longrightarrow} & H^{2}(H_{k})\oplus H^{4}(H_{k}) \\
\rotatebox{90}{$\in$} & & \rotatebox{90}{$\in$} \\[-4pt]
\xi & \longmapsto & c_{1}(\xi)\oplus c_{2}(\xi)
\end{array}
\end{eqnarray*}
is bijective.
\end{lemma}
\begin{proof}
Since $\dim_{\R} H_{k}=4$, any two bundles $\eta_{1}$ and $\eta_{2}\in {\rm Vect}_{\C}^{2}(H_{k})$ are isomorphic if and only if they are
stably isomorphic, i.e., $\eta_{1}\oplus \epsilon^\ell\equiv \eta_{2}\oplus \epsilon^\ell$ for some trivial complex $\ell$-dimensional bundle $\epsilon^{\ell}$,
see \cite[1.5 Theorem in Chapter 9]{Hu}.
Therefore $\eta_{1}$ and $\eta_{2}$ represent the same element in $\widetilde{K}(H_{k})$, the stable K-ring of $H_{k}$,  if and only if $\eta_{1}\equiv \eta_{2}$.
Therefore the map ${\rm Vect}_{\C}^{2}(H_{k})\to \widetilde{K}(H_{k})$ defined by $\xi\mapsto [\xi]$ is bijective.
Hence, it is enough to prove that the induced map
$$c': \widetilde{K}(H_{k})\to H^2(H_k)\oplus H^4(H_k),\qquad [\xi]\mapsto (c_1(\xi), c_2(\xi))$$
is bijective.

Let $s:\C P^1\to H_{k}=P(\gamma_{1}^k\oplus \epsilon^1)$ be the section defined by $s([p])=[p,[0:1]]$, and let $i:\C P^1\to H_{k}$ be an inclusion to a fiber in
the fibration $H_k\to \C P^1$.
Then $s(\C P^1)\cup i(\C P^1)\cong \C P^1\vee \C P^1$, and we have the following inclusion and collapsing sequence
\begin{eqnarray*}
\C P^1\vee \C P^1\longrightarrow H_{k}\longrightarrow H_{k}/(\C P^1\vee \C P^1).
\end{eqnarray*}
Since $H_k$ admits a CW-structure with one $0$-cell, two $2$-cells, and one $4$-cell (e.g.\ see \cite{DaJa}),
$H_{k}/(\C P^1\vee\C P^1)$ may be regarded as the collapsing of two $2$-cells to the one $0$-cell.
Therefore, the space $H_{k}/(\C P^1\vee\C P^1)$ is homeomorphic to $S^4$.
Hence, we have the following exact sequence of reduced $K$ groups (see \cite[2.1 Proposition in Chapter 10]{Hu}):
\begin{eqnarray*}
\widetilde{K}(S^4)\to \widetilde{K}(H_{k})\to \widetilde{K}(\C P^1\vee \C P^1).
\end{eqnarray*}

As is well known, we have the following isomorphisms 
\begin{eqnarray}
& &\widetilde{K}(S^4)\simeq \widetilde{K}(S^2)\simeq \widetilde{K}(\C P^1)\simeq \Z,\  {\rm and} \label{K-group of 4-sphere}\\
& &\widetilde{K}(\C P^1\vee \C P^1)\simeq \widetilde{K}(\C P^1)\oplus \widetilde{K}(\C P^1)\simeq \Z\oplus\Z=\Z^{2}.\label{K-group of 2-sphere}
\end{eqnarray}
These isomorphisms are induced by taking the Chern classes of vector bundles.
Let $c'=(c'_1, c_2'):\widetilde{K}(H_{k})\to H^{2}(H_{k})\oplus H^{4}(H_{k})\simeq  \Z^{2}\oplus\Z$, where $c'_1([\xi])=c_1(\xi)$ and $c'_2([\xi])=c_2(\xi)$.
Then $c'_1:\widetilde{K}(H_{k})\to H^{2}(H_{k})$ is surjective because for any $\alpha\in H^2(H_k)\simeq\Z^2$ can be realized as the first Chern class $c_1(\gamma)$ of a
complex line bundle $\gamma$ over $H_k$. 
Indeed, for a given $\alpha_{1}x+\alpha_{2}y\in \Z x\oplus\Z y=H^{2}(H_{k})$, the line bundle
$\gamma=\pi^{*}(\gamma_{1}^{\alpha_{1}})\otimes \gamma_{H_{k}}^{\alpha_{2}}$ has the first Chern class $\alpha_{1}x+\alpha_{2}y$, where
$\pi:H_{k}\to \C P^1$ is the projection, $\gamma_{H_{k}}$ is the canonical line bundle over $H_{k}=P(\gamma_{1}^{k}\oplus\epsilon^{1})$ induced from the vector bundle $\pi^{*}(\gamma_{1}^{k}\oplus\epsilon^{1})$, and $x$, $y$ are generators induced by $c_{1}(\pi^{*}\gamma_{1})$, $c_{1}(\gamma_{H_{k}})$ respectively.
We also claim that $c'_2:\widetilde{K}(H_{k})\to H^{4}(H_{k})$ is surjective. 
By the fundamental results of fibre bundle, 
we can construct all complex $2$-dimensional vector bundles over $H_{k}/(\C P^1\vee \C P^1)\cong S^{4}$
by using the continuous map $S^{4}\to BU(2)$ up to homotopy.
Because $\pi_{4}(BU(2))\simeq \Z$, for a given $\beta\in H^{4}(H_{k}/(\C P^1\vee \C P^1))$ we can construct the 
complex $2$-dimensional vector bundle $\eta'$ such that $c(\eta')=1+\beta$.
Now 
the collapsing map $\rho:H_{k}\to H_{k}/(\C P^1\vee \C P^1)$ induces the isomorphism 
$H^{4}(H_{k}/(\C P^1\vee \C P^1))\simeq H^{4}(H_{k})\simeq \Z$;
therefore, its pull-back $\eta=\rho^{*}\eta'$ over $H_{k}$ satisfies $c(\eta)=1+\beta$.
This implies that $c_{2}'$ is surjective.
Because $\gamma\oplus\eta$ is a complex $3$-dimensional vector bundle and $\dim_{\R}H_{k}=4$,
the bundle $\gamma\oplus\eta$ is in the stable range.
Therefore, there is the complex  $2$-dimensional vector bundle $\xi$ such that
$\xi\oplus\epsilon^1\equiv \gamma\oplus\eta$, where $\epsilon^1$ is the trivial line bundle over $H_{k}$,
and $c(\xi)=c(\gamma\oplus\eta)=1+c_{1}(\gamma)+c_{2}(\eta)$.
Therefore,
the map $c': \widetilde{K}(H_{k})\to H^2(H_k)\oplus H^4(H_k)$ is surjective.
Now consider the following diagram.
\begin{eqnarray*}
\begin{array}{ccccccccc}
  &                         & \widetilde{K}(S^4) & \longrightarrow & \widetilde{K}(H_{k}) & \longrightarrow & \widetilde{K}(\C P^1\vee \C P^1) & & \\
  &                         & \downarrow &  & \downarrow &  & \downarrow & & \\
0 & \longrightarrow & \Z                      & \longrightarrow & \Z^{2}\oplus\Z      & \longrightarrow & \Z^{2}                                         & \longrightarrow & 0
\end{array}
\end{eqnarray*}
Here the vertical maps from the left are  the isomorphism in \eqref{K-group of 4-sphere}, 
the map $c':\widetilde{K}(H_{k})\to H^2(H_k)\oplus H^4(H_k)$ and the isomorphism in \eqref{K-group of 2-sphere}, and the horizontal sequences are exact.
One can see easily that the diagram is commutative.
From the commutativity of the diagram and the surjectivity of the map $c'$, we can see that $\widetilde{K}(S^4) \to \widetilde{K}(H_{k}) \to \widetilde{K}(\C P^1\vee \C P^1)$ is a short exact sequence, and the map $c'$
is bijective.
Consequently, there exists the bijective map ${\rm Vect}^2_{\C}(H_{k})\to H^{2}(H_{k})\oplus H^{4}(H_{k})$ defined by $\xi\mapsto c_{1}(\xi)\oplus c_{2}(\xi)$.
This establishes the lemma.
\end{proof}

By Lemma \ref{basic classification}, any complex  $2$-dimensional  vector bundles over $H_{0}$ and $H_{1}$ can be written by
$$\eta_{(s,r,\alpha)}\to H_{0}, \ {\rm and } \quad\xi_{(s,r, \beta)}\to H_{1}$$
where
\begin{eqnarray*}
c_{1}(\eta_{(s,r,\alpha)})=(s,r)\in H^{2}(H_{0})\simeq \Z\oplus\Z,\quad c_{2}(\eta_{(s,r,\alpha)})=\alpha\in H^{4}(H_{0})\simeq \Z;\\
c_{1}(\xi_{(s,r,\beta)})=(s,r)\in H^{2}(H_{1})\simeq \Z\oplus\Z,\quad c_{2}(\xi_{(s,r,\beta)})=\beta\in H^{4}(H_{1})\simeq \Z.
\end{eqnarray*}
Moreover, by taking  tensor product with an appropriate line bundle if necessary, we may assume $(s,r)\in \{0,1\}^2$,  see Lemma \ref{formula of Chern class}.
Let $\mathcal{M}_{3}^{6}$ be the set of all $6$-dimensional $\C P$-towers of height $3$, up to diffeomorphism.
The main theorem of this section is the following.
\begin{theorem}
\label{3-stage}
The set $\mathcal{M}_{3}^{6}$ consists of the following distinct manifolds:
\begin{eqnarray*}
& &P(\eta_{(0,0,\alpha)})\ \text{for $\alpha\in \Z_{\ge 0}$}; \\
& &P(\eta_{(1,0,\alpha)})\ \text{for $\alpha\in \Z_{\ge 0}$}; \\
& &P(\eta_{(1,1,\alpha)})\ \text{for $\alpha\in \N$}; \\
& &P(\xi_{(0,0,\beta)})\ \text{for $\beta\in \N$}; \\
& &P(\xi_{(1,0,\beta)})\ \text{for $\beta\in \Z_{\ge 0}$}; \\
& &P(\xi_{(0,1,\beta)})\ \text{for $\beta\in \Z$}.
\end{eqnarray*}
Moreover, we have the  diffeomorphisms
$P(\eta_{(1,0,\alpha)})\cong P(\eta_{(0,1,\alpha)})$,
$P(\eta_{(0,0,1)})\cong P(\xi_{(0,0,0)})$, and
$P(\xi_{(0,1,\beta)})\cong P(\xi_{(1,1,-\beta)})$.
\end{theorem}
To prove Theorem \ref{3-stage}, we first observe the following.
For $H_{0}=\C P^1\times \C P^1$,
there is a self-diffeomorphism on $H_{0}$ defined by exchanging the first and second terms, i.e.,
$(p,q)\mapsto (q,p)$ for $(p,q)\in H_0= \C P^1\times \C P^1$.
This diffeomorphism induces a bundle isomorphism between $\eta_{(s,r,\alpha)}$ and $\eta_{(r,s,\alpha)}$.
Therefore, we may assume $(s,r)=(0,0),\ (1,0)$ or $(1,1)$ in the case of $\eta_{(s,r,\alpha)}$.

We also need the following lemma.
\begin{lemma}
\label{difference}
If  the cohomology ring $H^{*}(P(\eta_{(s,r,\alpha)}))$ is isomorphic to $H^{*}(P(\xi_{(s',r',\beta)}))$, then
$(s,r,\alpha)=(1,0,0)$ and $(s',r',\beta)=(0,0,0)$. Furthermore,  $P(\eta_{(1,0,0)})$ is diffeomorphic to $ P(\xi_{(0,0,0)})$.
\end{lemma}
\begin{proof}
By the Borel Hirzebruch formula \eqref{BH-formula}, we have the isomorphisms
\begin{eqnarray*}
H^{*}(P(\eta_{(s,r,\alpha)})) & \simeq & \Z[X,Y,Z]/\langle X^2,\ Y^2,\ Z^2+s ZX+r ZY+\alpha XY \rangle,\  \textrm{and}\\
H^{*}(P(\xi_{(s',r',\beta)})) & \simeq & \Z[x,y,z]/\langle x^2,\ y^2+xy,\ z^2+s' zx+r' zy+\beta xy \rangle,
\end{eqnarray*}
where $(s,r)=(0,0),\ (1,0)$ or $(1,1)$ in $\eta_{(s,r,\alpha)}$, and
$(s',r')=(0,0),\ (1,0),\ (0,1)$ or $(1,1)$ in $\xi_{(s',r',\alpha)}$.
For each $(s,r,\alpha)$ and $(s',r',\beta)$, we express the $\Z$-module structures of the above cohomology rings using their generators as follows:
\begin{eqnarray*}
& &\Z\oplus \Z X\oplus\Z Y\oplus \Z Z\oplus \Z XY \oplus \Z YZ \oplus \Z ZX \oplus \Z XYZ; \\
& &\Z\oplus \Z x\oplus\Z y\oplus \Z z\oplus \Z xy \oplus \Z yz \oplus \Z zx \oplus \Z xyz.
\end{eqnarray*}

Assume there exists an isomorphism $f:H^{*}(P(\eta_{(s,r,\alpha)}))\to H^{*}(P(\xi_{(s',r',\beta)}))$.
Let $f(X)=a_{1}x+b_{1}y+c_{1}z$, $f(Y)=a_{2}x+b_{2}y+c_{2}z$ and $f(Z)=a_{3}x+b_{3}y+c_{3}z$, and let $A_f$ denote
the corresponding $3\times 3$ matrix of $f$.
Because $f$ is a graded ring isomorphism, it satisfies the following relations:
\begin{eqnarray*}
f(X)^2=(a_{1}x+b_{1}y+c_{1}z)^2=(2a_{1}b_{1}-b_{1}^2-\beta c_{1}^{2})xy+(2a_{1}c_{1}-s'c_{1}^{2})xz+(2b_{1}c_{1}-r'c_{1}^{2})yz=0; \\
f(Y)^2=(a_{2}x+b_{2}y+c_{2}z)^2=(2a_{2}b_{2}-b_{2}^2-\beta c_{2}^{2})xy+(2a_{2}c_{2}-s'c_{2}^{2})xz+(2b_{2}c_{2}-r'c_{2}^{2})yz=0
\end{eqnarray*}
in $H^{*}(P(\xi_{(s',r',\beta)}))$. Therefore, we have 
\begin{eqnarray*}
& &2a_{i}b_{i}-b_{i}^2-\beta c_{i}^{2}=0; \\
& &2a_{i}c_{i}-s'c_{i}^{2}=0; \\
& &2b_{i}c_{i}-r'c_{i}^{2}=0,
\end{eqnarray*}
for $i=1,2$.

Assume $c_{1}=0$.
Then, by using the first equation above and $\det A_{f}=\pm 1$, we have either $b_{1}=0$ and $\ a_{1}=\epsilon_{1}$, or $b_{1}=2a_{1}=2\epsilon_{1}$, where $\epsilon_{1}=\pm 1$.
If $c_{2}=0$, then it is easy to check that this gives a contradiction to $\det A_{f}=\pm 1$.
Hence, $c_{2}\not=0$.
By using the second and the third equations above, we have  $s'c_{2}=2a_{2}$ and $r'c_{2}=2b_{2}$.
Hence it can be seen easily from $\det A_{f}=\pm 1$ that only $(s',r')=(0,0)$ is possible,
and in this case $(a_{2},b_{2},c_{2})=(0,0,\epsilon_{2})$ and $\beta=0$, where $\epsilon_{2}=\pm 1$.
Hence, we have that $(s',r',\beta)=(0,0,0)$.

If $(a_{1},b_{1},c_{1})=(\epsilon_{1},0,0)$, then $b_{3}=\epsilon_{3}$ because $\det A_{f}=\pm 1$. Therefore, it follows from $f(Z)^2=-s f(X)f(Z)-r f(Y)f(Z)-\alpha f(X)f(Y)$ that
\begin{eqnarray*}
& &2a_{3}\epsilon_{3}-1=-s\epsilon_{1}\epsilon_{3}; \\
& &2a_{3}c_{3}=-s\epsilon_{1}c_{3}-r\epsilon_{2}a_{3}-\alpha\epsilon_{1}\epsilon_{2}; \\
& &2\epsilon_{3}c_{3}=-r\epsilon_{2}\epsilon_{3}.
\end{eqnarray*}
Using the third equation above, we have  $r=c_{3}=0$.
Therefore, by the second equation, we also have $\alpha=0$.
Moreover, from the first equation $s=1$.
Hence, $(s,r,\alpha)=(1,0,0)$.

If $(a_{1},b_{1},c_{1})=(\epsilon_{1},2\epsilon_{1},0)$, then $b_{3}-2a_{3}=\epsilon_{3}$ because $\det A_{f}=\pm 1$.
Therefore, it follows from $f(Z)^2=-s f(X)f(Z)-r f(Y)f(Z)-\alpha f(X)f(Y)$ that
\begin{eqnarray*}
& &2a_{3}b_{3}-b_{3}^2=s\epsilon_{1}b_{3}-2s\epsilon_{1}a_{3}; \\
& &2a_{3}c_{3}=-s\epsilon_{1}c_{3}-r\epsilon_{2}a_{3}-\alpha\epsilon_{1}\epsilon_{2}; \\
& &2b_{3}c_{3}=-r\epsilon_{2}b_{3}-2s\epsilon_{1}c_{3}-2\alpha \epsilon_{1}\epsilon_{2}.
\end{eqnarray*}
Using the first equation and $b_{3}-2a_{3}=\epsilon_{3}$, we have $b_{3}=-s\epsilon_{1}$.
Therefore, by using the third equation, we have  $sr=-2\alpha$.
This implies that $\alpha=0$ and $sr=0$.
If $s=0$, then $b_{3}=-s\epsilon_{1}=0$; however, $b_{3}-2a_{3}=-2a_{3}=\epsilon_{3}$ and this gives a contradiction.
Therefore $(s,r,\alpha)=(1,0,0)$.
This establishes the first statement of the lemma when $c_1=0$ case.

In the case when $c_{1}\not=0$ and  $c_{2}=0$,
by a similar argument to the above case, we have the same result.
When $c_{1}\not=0$ and $c_{2}\not=0$,
by some routine computation, we can see that this case gives a contradiction.
This establishes the first statement of the lemma.

Because $\eta_{(1,0,0)}\equiv \gamma_{x}\oplus \epsilon$, where $\gamma_{x}$ is the tautological line bundle along the first factor of $\C P^1\times \C P^1$, we can easily check that $P(\eta_{(1,0,0)})\cong (S^{3}\times \C P^1) \times_{T^1}P(\C_{1}\oplus\C)$,
where $T^1$ acts on $S^3$ as diagonal multiplications in its coordinates and  trivially on $\C P^1$ and $\C_1$ is a complex $1$-dimensional $T^1$ representation such that $t\cdot z=tz$ for $t\in T^1$ and $z\in \C_1$.
On the other hand, because $\xi_{(0,0,0)}$ is the trivial bundle over $H_{1}$ (by Lemma \ref{basic classification}), we have that$P(\xi_{(0,0,0)})=S^{3}\times_{T^1}P(\C_{1}\oplus\C)\times \C P^1$.
Therefore, we have that $P(\eta_{(1,0,0)})\cong P(\xi_{(0,0,0)})$.
This establishes the second statement.
\end{proof}

In order to prove Theorem \ref{3-stage}, we may divide the proof into the following two cases.

\begin{description}
\item[CASE I] $P(\eta_{(s,r,\alpha)})$ with the base space $H_0$. In this case  $(s,r)=(0,0), (1,0)$ and $(1,1)$.
\item[CASE II] $P(\xi_{(s,r,\alpha)})$ with the base space $H_1$. In this case $(s,r)=(0,0), (1,0), (0,1)$ and $(1,1)$. Moreover if $(s,r)=(0,0)$ then $\alpha\not=0$.
\end{description}

The rest of the section in devoted to the proof of Theorem \ref{3-stage} by treating the two cases separately.

\subsection*{CASE I: $P(\eta_{(s,r,\alpha)})$ with the base space $H_{0}$}

We prove the cohomological rigidity for $P(\eta_{(s,r,\alpha)})$. 
Namely, we prove the following proposition.
\begin{proposition}
\label{cohomological rigidity 1}
The following statements are equivalent.
\begin{enumerate}
\item Two manifolds  $P(\eta_{(s_{1},r_{1},\alpha_{1})})$ and $ P(\eta_{(s_{2},r_{2},\alpha_{2})})$ are diffeomorphic.
\item Two cohomology rings $H^{*}(P(\eta_{(s_{1},r_{1},\alpha_{1})}))$ and  $H^{*}(P(\eta_{(s_{2},r_{2},\alpha_{2})}))$ are isomorphic. 
\item $(s_{1},r_{1})=(s_{2},r_{2})$,  and $\alpha_1$ and $\alpha_2$ are as follows:
\begin{enumerate}
\item if   $(s_{1},r_{1})=(s_{2},r_{2})=(0,0)$, then $\alpha_{2}=\alpha_{1}$ or $-\alpha_{1}$;
\item if $(s_{1},r_{1})=(s_{2},r_{2})=(1,0)$ (or $(0,1)$), then $\alpha_{2}=\alpha_{1}$ or $-\alpha_{1}$;
\item if $(s_{1},r_{1})=(s_{2},r_{2})=(1,1)$, then  $\alpha_{2}=\alpha_{1}$ or $-\alpha_{1}+1$.
\end{enumerate}
\end{enumerate}
\end{proposition}
\begin{proof}
$(1)\Rightarrow (2)$ is trivial.

We first prove $(2)\Rightarrow (3)$.
By  \eqref{BH-formula}, we have the following isomorphisms 
\begin{eqnarray*}
H^{*}(P(\eta_{(s_1,r_1,\alpha_1)})) & \simeq & \Z[X,Y,Z]/\langle X^2,\ Y^2,\ Z^2+s_1 ZX+r_1 ZY+\alpha_{1} XY \rangle,\  \textrm{and}\\
H^{*}(P(\eta_{(s_2,r_2,\alpha_2)}))& \simeq & \Z[x,y,z]/\langle x^2,\ y^2,\ z^2+s_2 zx+r_2 zy+\alpha_2xy \rangle.
\end{eqnarray*}
Assume there exists a graded ring isomorphism $f:H^{*}(P(\eta_{(s_{1},r_{1},\alpha_{1})}))\simeq H^{*}(P(\eta_{(s_{2},r_{2},\alpha_{2})}))$,
and put
the matrix representation of $f:H^{2}(P(\eta_{(s_{1},r_{1},\alpha_{1})}))\simeq H^{2}(P(\eta_{(s_{2},r_{2},\alpha_{2})}))$ with respect to the given module generators as
\begin{eqnarray*}
A_{f}=
\left(
\begin{array}{ccc}
a_{1} & b_{1} & c_{1} \\
a_{2} & b_{2} & c_{2} \\
a_{3} & b_{3} & c_{3} \\
\end{array}
\right),
\end{eqnarray*}
i.e.,
$f(X)=a_{1}x+b_{1}y+c_{1}z$, $f(Y)=a_{2}x+b_{2}y+c_{2}z$, $f(Z)=a_{3}x+b_{3}y+c_{3}z$.
Note that $\det A_{f}=\pm 1$.
Because $X^2=Y^2=0$ and $f$ is a ring isomorphism,
\begin{eqnarray*}
& &f(X)^2=(2a_{1}b_{1}-\alpha_{2}c_{1}^2)xy+(2a_{1}-s_{2}c_{1})c_{1}xz+(2b_{1}-r_{2}c_{1})c_{1}yz=0; \\
& &f(Y)^2=(2a_{2}b_{2}-\alpha_{2}c_{2}^2)xy+(2a_{2}-s_{2}c_{2})c_{2}xz+(2b_{2}-r_{2}c_{2})c_{2}yz=0
\end{eqnarray*}
in $H^{*}(P(\eta_{(s_2,r_2,\alpha_2)}))$. Therefore, we have
\begin{eqnarray}
& &\label{eq-1}2a_{i}b_{i}-\alpha_{2}c_{i}^2=0, \\
& &\label{eq-2}(2a_{i}-s_{2}c_{i})c_{i}=0, \\
& &\label{eq-3}(2b_{i}-r_{2}c_{i})c_{i}=0,
\end{eqnarray}
for $i=1,\ 2$.
We divide the proof into the following three cases: 
{\bf Case 1} $(s_{2},r_{2})=(1,1)$; {\bf Case 2} $(s_{2},r_{2})=(0,0)$; {\bf Case 3}  $(s_{2},r_{2})=(1,0)$.

\textbf{Case 1: $(s_{2},r_{2})=(1,1)$.}
We first claim that $c_{1}=c_{2}=0$ and $c_{3}=\epsilon_{3}=\pm 1$.
If $c_{i}\not=0$, for $i=1$ or $2$, then $2a_{i}=c_{i}$ by \eqref{eq-2}, $2b_{i}=c_{i}$ by \eqref{eq-3} and
$2a_{i}b_{i}=\alpha_{2}c_{i}^2$ by \eqref{eq-1}.
These equations imply that
\begin{eqnarray*}
4a_{i}b_{i}=c_{i}^2=2\alpha_{2}c_{i}^2.
\end{eqnarray*}
Because $c_{i}\not=0$, we have that $1=2\alpha_{2}$.
This gives a contradiction.
Therefore, we have 
\begin{eqnarray*}
c_{1}=c_{2}=0.
\end{eqnarray*}
This together with $\det A_{f}=\pm 1$ imply that  
\begin{eqnarray*}
c_{3}=\epsilon_{3}=\pm 1.
\end{eqnarray*}

Because $Z^2=-s_{1}XZ-r_{1}YZ-\alpha_{1}XY$ in $H^{*}(P(\eta_{(s_{1},r_{1},\alpha_{1})}))$,
the ring isomorphism $f$ induces the following equations
\begin{eqnarray}
& &\label{eq-4}2a_{3}b_{3}-\alpha_{2}\epsilon_{3}^2=-s_{1}(a_{1}b_{3}+a_{3}b_{1})-r_{1}(a_{2}b_{3}+a_{3}b_{2})-\alpha_{1}(a_{1}b_{2}+a_{2}b_{1}),\\
& &\label{eq-5}(2a_{3}-\epsilon_{3})\epsilon_{3}=(-s_{1}a_{1}-r_{1}a_{2})\epsilon_{3},\\
& &\label{eq-6}(2b_{3}-\epsilon_{3})\epsilon_{3}=(-s_{1}b_{1}-r_{1}b_{2})\epsilon_{3}.
\end{eqnarray}
Using \eqref{eq-1} and $c_{1}=c_{2}=0$,
we have $a_{i}b_{i}=0$ for $i=1,\ 2$.
Moreover, from $\det A_{f}=\pm 1$,
there are two possibilities, i.e.,  either $(a_{1},b_{2})=(0,0)$ and $(a_{2},b_{1})=(\epsilon_{1},\epsilon_{2})$, or  $(a_{1},b_{2})=(\epsilon_{1},\epsilon_{2})$ and $(a_{2},b_{1})=(0,0)$ where $\epsilon_{i}=\pm 1$ for $i=1,2$.

If $(a_{1},b_{2})=(0,0)$ and $(a_{2},b_{1})=(\epsilon_{1},\epsilon_{2})$, then
it follows from \eqref{eq-5} and \eqref{eq-6} that
\begin{eqnarray*}
& &2a_{3}=\epsilon_{3}-r_{1}\epsilon_{1}; \\
& &2b_{3}=\epsilon_{3}-s_{1}\epsilon_{2}.
\end{eqnarray*}
It is easy to check that if $s_{1}=0$ or $r_{1}=0$ then we have a contradiction to one of the equations above.
Therefore, $(s_{1},r_{1})=(s_{2},r_{2})=(1,1)$.
We also have that
if $\epsilon_{3}=\epsilon_{1}$ (resp.\ $\epsilon_{3}=\epsilon_{2}$) then $a_{3}=0$ (resp.\ $b_{3}=0$)
and if $\epsilon_{3}\not=\epsilon_{1}$ (resp.\ $\epsilon_{3}\not=\epsilon_{2}$) then $a_{3}=\epsilon_{3}$ (resp.\ $b_{3}=\epsilon_{3}$).
Thus, by the equation \eqref{eq-4}, we have that $\alpha_{2}=\alpha_{1}$ or $\alpha_{2}=-\alpha_{1}+1$.

If $(a_{1},b_{2})=(\epsilon_{1},\epsilon_{2})$ and $(a_{2},b_{1})=(0,0)$, then similarly we have that
$(s_{1},r_{1})=(s_{2},r_{2})=(1,1)$ and $\alpha_{2}=\alpha_{1}$ or $\alpha_{2}=-\alpha_{1}+1$.
This establishes $(3)-(c)$.

\textbf{Case 2: $(s_{2},r_{2})=(0,0)$.}
If $(s_1, r_1)=(1,1)$ in this case, by the same argument as in Case 1 with $(s_2, r_2)$ replaced by $(s_1, r_1)$, we can see that
$(s_2, r_2)=(1,1)$ which contradicts to the hypothesis. Therefore 
$(s_{1},r_{1})=(0,0)$ or $(1,0)$, and hence, $Z^2=-s_{1}XZ-\alpha_{1}XY$ in $H^{*}(P(\eta_{(s_{1},r_{1},\alpha_{1})}))$.
Therefore, the ring isomorphism $f$ implies the following equations:
\begin{eqnarray}
& &\label{eq-4'}2a_{3}b_{3}-\alpha_{2}c_{3}^2=-s_{1}(a_{1}b_{3}+a_{3}b_{1})-\alpha_{1}(a_{1}b_{2}+a_{2}b_{1})+s_{1}c_{1}c_{3}\alpha_{2}+\alpha_{1}c_{1}c_{2}\alpha_{2}; \\
& &\label{eq-5'}2a_{3}c_{3}=-s_{1}(a_{1}c_{3}+a_{3}c_{1})-\alpha_{1}(a_{1}c_{2}+a_{2}c_{1}); \\
& &\label{eq-6'}2b_{3}c_{3}=-s_{1}(b_{1}c_{3}+b_{3}c_{1})-\alpha_{1}(b_{1}c_{2}+b_{2}c_{1}).
\end{eqnarray}

Because of \eqref{eq-2} and \eqref{eq-3}, we also have that $a_{i}c_{i}=b_{i}c_{i}=0$.
Then by \eqref{eq-1}, there are  two cases to consider for $i=1,2$:
\textbf{(2-i)} the case when  $c_{i}\not=0$, and hence, $a_{i}=b_{i}=\alpha_{2}=0$; 
\textbf{(2-ii)} the case when   $c_{i}=0$, and hence $a_{i}b_{i}=0$.

\textbf{(2-i)} If $c_{1}\not=0$, and hence, $a_{1}=b_{1}=\alpha_{2}=0$, then $c_{1}=\epsilon_{3}=\pm 1$ because $\det A_{f}=\pm 1$.
Furthermore, if $c_{2}\not=0$, then $a_{2}=b_{2}=0$,  which gives a contradiction to  $\det A_{f}=\pm 1$.
Therefore, $c_{2}=0$ and $a_{2}b_{2}=0$. Moreover  $a_{3}b_{3}=0$ by  \eqref{eq-4'}.
Since  $\det A_{f}=\pm 1$, there are two possibilities for $(a_{2},a_{3})$ and  $(b_{2},b_{3})$, i.e.,  either 
$(a_{2},a_{3})=(0,\epsilon_{1})$ and $(b_{2},b_{3})=(\epsilon_{2},0)$,  or 
$(a_{2},a_{3})=(\epsilon_{1},0)$ and $(b_{2},b_{3})=(0,\epsilon_{2})$.
If $a_{2}=b_{3}=0$, then, by using \eqref{eq-5'} and \eqref{eq-6'}, we have that
$2c_{3}=-s_{1}\epsilon_{3}$ and $\alpha_{1}=\alpha_{2}=0$.
Therefore, because $s_{1}=0$ or $1$, we also have $c_{3}=0$ and $s_{1}=s_{2}=0$.
If $a_{3}=b_{2}=0$, then we similarly have that $\alpha_{1}=\alpha_{2}=0$ and $s_{1}=s_{2}=0$.

\textbf{(2-ii)}
If $c_{1}=0$, then $a_{1}b_{1}=0$.
If $c_{2}\not=0$, then the proof is almost the same with the case when $c_{1}\not=0$; and we have that $\alpha_{1}=\alpha_{2}=0$ and $s_{1}=s_{2}=0$ as the conclusion.
Therefore, we may put $c_{2}=0$ and $a_{2}b_{2}=0$.
Because of $\det A_{f}=\pm 1$, we have that $c_{3}=\epsilon_{3}=\pm 1$ and there are the  two possibilities, i.e., either 
$(a_{1},a_{2})=(0,\epsilon_{1})$ and $(b_{1},b_{2})=(\epsilon_{2},0)$, or $(a_{1},a_{2})=(\epsilon_{1},0)$ and $(b_{1},b_{2})=(0,\epsilon_{2})$.
If $a_{1}=b_{2}=0$ (resp.\ $a_{2}=b_{1}=0$), then it follows from \eqref{eq-6'} (resp.\ \eqref{eq-5'}) that $2b_{3}=-s_{1}b_{1}$ (resp.\ $2a_{3}=-s_{1}a_{1}$).
Therefore, $s_{1}=s_{2}=0$ and $b_{3}=0$ (resp.\ $a_{3}=0$).
Moreover, by \eqref{eq-4'}, we have that $\alpha_{2}=\epsilon_{1}\epsilon_{2}\alpha_{1}$.
This establishes $(3)-(a)$.

\textbf{Case 3: $(s_{2},r_{2})=(1,0)$.}
In this case, by the same arguments as above, we may assume $(s_{1},r_{1})=(1,0)$, i.e., $Z^2=-XZ-\alpha_{1}XY$ in $H^{*}(P(\eta_{(s_{1},r_{1},\alpha_{1})}))$.
It is sufficient to show that $\alpha_{2}=\alpha_{1}$ or $-\alpha_{1}$.
Now, the ring isomorphism $f$ implies the following equations:
\begin{eqnarray}
& &\label{eq-4''}2a_{3}b_{3}-\alpha_{2}c_{3}^2=-(a_{1}b_{3}+a_{3}b_{1})-\alpha_{1}(a_{1}b_{2}+a_{2}b_{1})+c_{1}c_{3}\alpha_{2}+\alpha_{1}c_{1}c_{2}\alpha_{2}; \\
& &\label{eq-5''}2a_{3}c_{3}-c_{3}^2=-(a_{1}c_{3}+a_{3}c_{1})-\alpha_{1}(a_{1}c_{2}+a_{2}c_{1})+c_{1}c_{3}+c_{1}c_{2}\alpha_{1}; \\
& &\label{eq-6''}2b_{3}c_{3}=-(b_{1}c_{3}+b_{3}c_{1})-\alpha_{1}(b_{1}c_{2}+b_{2}c_{1}).
\end{eqnarray}
Because of \eqref{eq-2} and \eqref{eq-3}, we also have $(2a_{i}-c_{i})c_{i}=0$ and $b_{i}c_{i}=0$.
By \eqref{eq-1}, if $c_{i}\not=0$ then $b_{i}=\alpha_{2}=0$ and  $c_{i}=2a_{i}$, and if $c_{i}=0$, then $a_{i}b_{i}=0$.

\textbf{(3-i)}
If $c_{1}\not=0$, then $b_{1}=\alpha_{2}=0$, $c_{1}=2a_{1}$.
Since $\det A_{f}=\pm 1$, we may put $a_{1}=\epsilon_{1}=\pm 1$.
In this case, if $c_{2}\not=0$ then $b_{2}=0$ and $c_{2}=2a_{2}$, which contradicts to $\det A_{f}=\pm 1$.
Therefore, $c_{2}=0$ and $a_{2}b_{2}=0$.
It follows from \eqref{eq-4''} and \eqref{eq-6''} that
\begin{eqnarray*}
2a_{3}b_{3}=-\epsilon_{1}(b_{3}+\alpha_{1}b_{2})=b_{3}c_{3}.
\end{eqnarray*}
Therefore, there are  two cases to consider: the case when $b_{3}=0$,  and  hence $\alpha_{1}b_{2}=0$; the case when
$b_{3}\not=0$,  and hence $c_{3}=2a_{3}$.
If  $b_{3}\not=0$ and $c_{3}=2a_{3}$, then 
by  $\det A_{f}=\pm 1$ we have $a_{3}=0=c_{3}$ and $b_{3}=\epsilon_{2}=\pm 1$.
Then the matrix $A_{f}$ is equal 
\begin{eqnarray*}
\left(
\begin{array}{ccc}
\epsilon_{1} & 0 & 2\epsilon_{1} \\
a_{2} & b_{2} & 0 \\
0 & \epsilon_{2} & 0
\end{array}
\right).
\end{eqnarray*}
This gives a contradiction to  $\det A_{f}=\pm 1$.
Therefore, $b_{3}=0$, and  hence $\alpha_{1}b_{2}=0$.
If $b_{2}=0$ then this gives a contradiction to $\det A_{f}=\pm 1$.
Hence, we have $b_{2}\not=0$, and hence $\alpha_{1}=\alpha_{2}=0$.

\textbf{(3-ii)} If $c_{1}=0$ and $c_{2}\not=0$, then 
$a_{1}b_{1}=0$, $c_{2}=2a_{2}$ and $b_{2}=\alpha_{2}=0$.
If $b_{1}=0$, then it is easy to check this gives a contradiction to $\det A_{f}=\pm 1$.
Hence, $a_{1}=0$ and $b_{1}=\pm 1$.
Because $c_{2}=2a_{2}$ and $\det A_{f}=\pm 1$, we have $c_{3}-2a_{3}=\pm 1$.
By using \eqref{eq-5''}, we also have the equation $c_{3}(c_{3}-2a_{3})=0$.
Therefore, $c_{3}=0$, and hence $2a_{3}=\pm 1$.
This gives a contradiction to $a_{3}\in \Z$.

Therefore $c_{1}=c_{2}=0$.
Since $\det A_{f}=\pm 1$ and $c_{1}=c_{2}=0$, we can put $c_{3}=\epsilon_{3}=\pm 1$.
Then, we can easily see that $a_{1}+2a_{3}=\epsilon_{3}$ by \eqref{eq-5''} and $b_{1}=-2b_{3}$ by \eqref{eq-6''}.
Therefore, by using $a_{1}b_{1}=a_{2}b_{2}=0$ and $\det A_{f}=\pm 1$,
we have that $b_{1}=b_{3}=0$, $b_{2}=\epsilon_{2}=\pm 1$ and $a_{2}=0$, $a_{1}=\epsilon_{1}=\pm 1$.
Hence, by using \eqref{eq-4''}, we have  $\alpha_{2}=\pm \alpha_{1}$.
This establishes $(3)-(b)$.
Consequently, we have proved the implication $(2)\Rightarrow (3)$.

Finally, we prove $(3)\Rightarrow (1)$.
Consider the diffeomorphism $f={\rm id}\times {\rm conj}:\C P^1\times \C P^1\to \C P^1\times \C P^1$ defined by $(p,q)\mapsto (p, \overline{q})$.
Because $f$ changes the orientation on $\C P^1\times \C P^1$, the Euler class  $e(f^{*}\eta_{(s,r,\alpha)})$ coincides with
$-e(\eta_{(s,r,\alpha)})$.
Because of the definition of Chern class, $e(f^{*}\eta_{(s,r,\alpha)})=c_{2}(f^{*}\eta_{(s,r,\alpha)})=-c_{2}(\eta_{(s,r,\alpha)})=-\alpha$.
Because $x$ and $y$ are the first Chern classes of the tautological line bundles of the first and the second  factor of $\C P^1\times \C P^1$, we have $c_{1}(f^{*}\eta_{(s,r,\alpha)})=f^{*}(sX+rY)=sx-ry$.
Hence, by Lemmas \ref{formula of Chern class} and \ref{basic classification}, we have
\begin{eqnarray*}
& &f^{*}\eta_{(s,0,\alpha)}\equiv \eta_{(s,0,-\alpha)}; \\
& &f^{*}\eta_{(1,1,\alpha)}\otimes \gamma_{2}\equiv \eta_{(1,-1,-\alpha)}\otimes \gamma_{2}\equiv  \eta_{(1,1,1-\alpha)},
\end{eqnarray*}
where $\gamma_{2}$ is the pull back of the tautological line bundle over $\C P^1$  along the projection $\pi_{2}:\C P^1\times \C P^1\to \C P^1$  to the second factor.
This implies that $P(\eta_{(s,r,\alpha)})\cong P(\eta_{(s,r,-\alpha)})$ for $(s,r)=(0,0)$ or $(1,0)$ (or $(0,1)$) and
$P(\eta_{(1,1,\alpha)})\cong P(\eta_{(1,1,1-\alpha)})$ for $(s,r)=(1,1)$.
This proves the implication $(3)\Rightarrow (1)$.
\end{proof}

\subsection*{CASE II: $P(\xi_{(s,r,\beta)})$ with the base space $H_{1}$}
We  prove the cohomological rigidity for
$P(\xi_{(s,r,\beta)})$ 
in the following proposition.
\begin{proposition}
\label{cohomological rigidity 2}
The following statements are equivalent.
\begin{enumerate}
\item Two manifolds $P(\xi_{(s_{1},r_{1},\beta_{1})})$ and $P(\xi_{(s_{2},r_{2},\beta_{2})})$  are  diffeomorphic.
\item Two cohomology rings $H^{*}(P(\xi_{(s_{1},r_{1},\beta_{1})}))$ and $ H^{*}(P(\xi_{(s_{2},r_{2},\beta_{2})}))$  are isomorphic.
\item  Either $(s_{1},r_{1},\beta_{1})=(s_{2},r_{2},\beta_{2})$, or one of the following holds:
\begin{enumerate}
\item $(s_{1},r_{1},\beta_{1})=(0,0,\beta)$ and $(s_{2},r_{2},\beta_{2})=(0,0,-\beta)$ ($\beta\not=0$);
\item $(s_{1},r_{1},\beta_{1})=(1,0,\beta)$ and $(s_{2},r_{2},\beta_{2})=(1,0,-\beta)$;
\item $\{(s_{1},r_{1},\beta_{1}),\ (s_{2},r_{2},\beta_{2})\}=\{(0,1,\beta),\ (1,1,-\beta)\}$,
\end{enumerate}
for some $\beta\in \Z$.
\end{enumerate}
\end{proposition}

By using Proposition \ref{cohomological rigidity 1} and \ref{cohomological rigidity 2} and Lemma \ref{difference},
we have Theorem \ref{3-stage}.
Let us prove Proposition \ref{cohomological rigidity 2}.
\begin{proof}
$(1)\Rightarrow (2)$ is trivial.
We first prove $(2)\Rightarrow (3)$.
By  \eqref{BH-formula} we have the isomorphisms
\begin{eqnarray*}
H^{*}(P(\xi_{(s_1,r_1,\beta_1)})) & \simeq & \Z[X,Y,Z]/\langle X^2,\ Y^2+XY,\ Z^2+s_1 ZX+r_1 ZY+\beta_1 XY \rangle,\  \textrm{and}\\
H^{*}(P(\xi_{(s_2,r_2,\beta_2)})) & \simeq & \Z[x,y,z]/\langle x^2,\ y^2+xy,\ z^2+s_2 zx+r_2 zy+\beta_2 xy \rangle.
\end{eqnarray*}
Assume  there is a ring isomorphism $f:H^{*}(P(\xi_{(s_{1},r_{1},\beta_{1})}))\simeq H^{*}(P(\xi_{(s_{2},r_{2},\beta_{2})}))$, and put
the matrix representation of $f:H^{2}(P(\xi_{(s_{1},r_{1},\beta_{1})}))\simeq H^{2}(P(\xi_{(s_{2},r_{2},\beta_{2})}))$ as
\begin{eqnarray*}
A_{f}=
\left(
\begin{array}{ccc}
a_{1} & b_{1} & c_{1} \\
a_{2} & b_{2} & c_{2} \\
a_{3} & b_{3} & c_{3} \\
\end{array}
\right).
\end{eqnarray*}
Note that $\det A_{f}=\pm 1$.
Let $\epsilon_{i}=\pm 1$ ($i=1,\ 2,\ 3$).
Because of $X^2=0 \in H^{*}(P(\xi_{(s_1,r_1,\beta_1)}))  $, we have 
\begin{eqnarray*}
& &2a_{1}b_{1}-b_{1}^2-c_{1}^2\beta_{2}=0, \\
& &2a_{1}c_{1}-c_{1}^2s_{2}=0, \\
& &2b_{1}c_{1}-c_{1}^2r_{2}=0.
\end{eqnarray*}
By using these equations and $\det A_{f}=\pm 1$, it is easy to check that for $\epsilon=\pm 1$
\begin{description}
\item[Case 1] if $c_{1}\not=0$, then there are the following two sub-cases:
\begin{itemize}
\item $(s_{2},r_{2})=(0,0)$ with $(a_{1},b_{1},c_{1})=(0,0,\epsilon)$ and $\beta_{2}=0$;
\item $(s_{2},r_{2})=(1,0)$ with $(a_{1},b_{1},c_{1})=(\epsilon,0,2\epsilon)$ and $\beta_{2}=0$,
\end{itemize}
\item[Case 2] if $c_{1}=0$, then $(a_{1},b_{1})=(\epsilon,0)$ or $(\epsilon,2\epsilon)$.
\end{description}
Because  $Y^2=-XY$ in $ H^{*}(P(\xi_{(s_1,r_1,\beta_1)})) $, we also have
\begin{eqnarray}
& &\label{eq2-4}2a_{2}b_{2}-b_{2}^2-c_{2}^2\beta_{2}=-a_{1}b_{2}-b_{1}a_{2}+b_{1}b_{2}+c_{1}c_{2}\beta_{2},\\
& &\label{eq2-5}2a_{2}c_{2}-c_{2}^2s_{2}=-a_{1}c_{2}-c_{1}a_{2}+c_{1}c_{2}s_{2},\\
& &\label{eq2-6}2b_{2}c_{2}-c_{2}^2r_{2}=-b_{1}c_{2}-c_{1}b_{2}+c_{1}c_{2}r_{2}.
\end{eqnarray}

\textbf{Case 1: $c_{1}\not=0$.}
If $(s_{2},r_{2})=(0,0)$, then, by using \eqref{eq2-5}, \eqref{eq2-6} and $(a_{1},b_{1},c_{1})=(0,0,\epsilon_{3})$,
we can easily show that $a_{2}=b_{2}=0$;
however, because $\det A_{f}=\pm 1$, this gives a contradiction.
Therefore, $(s_{2},r_{2},\beta_{2})=(1,0,0)$ and $(a_{1},b_{1},c_{1})=(\epsilon_{1},0,2\epsilon_{1})$.
Note that $\det A_{f} (a_{2}b_{3}-a_{3}b_{2})$ is the $(1,3)$-entry of  the matrix $A_{f}^{-1}$.
Therefore, by a similar argument to the above, we can see that  if $a_{2}b_{3}-a_{3}b_{2}\not=0$ then $(s_{1},r_{1})=(1,0)$ and $\beta_{1}=0$.
This means that if we get $a_{2}b_{3}-a_{3}b_{2}\not=0$ then we have $(s_{1},r_{1},\beta_{1})=(s_{2},r_{2},\beta_{2})=(1,0,0)$, i.e., 
the statement of this proposition holds.

By \eqref{eq2-6},  we may divide the case when $c_{1}\not=0$ into  two sub-cases: {\bf (1-i)} $b_{2}=0$ and  {\bf (1-ii)}   $b_{2}\not=0$ and $c_{2}=-\epsilon_{1}$. 

{\bf (1-i)}  If $b_{2}=0$, then it easily follows from \eqref{eq2-5} that $c_{2}=2a_{2}$ or $-\epsilon_{1}$.
Moreover, by using $\det A_{f}=\pm 1$ and $(a_{1},b_{1},c_{1})=(\epsilon_{1},0,2\epsilon_{1})$, we have that
$(a_{2},b_{2},c_{2})=(0,0,-\epsilon_{1})$ or $(-\epsilon_{1},0,-\epsilon_{1})$, and $b_{3}=\epsilon_{2}$.
If  $(a_{2},b_{2},c_{2})=(-\epsilon_{1},0,-\epsilon_{1})$, then $a_{2}b_{3}-a_{3}b_{2}=-\epsilon_{1}\epsilon_{2}\not=0$.
Therefore, by the argument explained above, we have  $(s_{1},r_{1},\beta_{1})=(s_{2},r_{2},\beta_{2})=(1,0,0)$.
Hence, this satisfies the statement of this proposition.
Suppose $(a_{2},b_{2},c_{2})=(0,0,-\epsilon_{1})$.
Since $Z^{2}=-s_{1}XZ-r_{1}YZ-\beta_{1}XY$, we have 
\begin{eqnarray*}
& &(2a_{3}\epsilon_{2}-1)xy+2\epsilon_{2}c_{3}yz+(2a_{3}c_{3}-c_{3}^2)xz \\
&=&-s_{1}(\epsilon_{1}x+2\epsilon_{1}z)(a_{3}x+\epsilon_{2}y+c_{3}z)
+r_{1}\epsilon_{1}z(a_{3}x+\epsilon_{2}y+c_{3}z)+\beta_{1}(\epsilon_{1}x+2\epsilon_{1}z)\epsilon_{1}z.
\end{eqnarray*}
So, we have
\begin{eqnarray*}
& &2a_{3}\epsilon_{2}-1=-s_{1}\epsilon_{1}\epsilon_{2}; \\
& &2a_{3}c_{3}-c_{3}^2=-2s_{1}\epsilon_{1}a_{3}+s_{1}\epsilon_{1}c_{3}+r_{1}\epsilon_{1}a_{3}-r_{1}\epsilon_{1}c_{3}-\beta_{1}; \\
& &2\epsilon_{2}c_{3}=-2s_{1}\epsilon_{1}\epsilon_{2}+r_{1}\epsilon_{1}\epsilon_{2}.
\end{eqnarray*}
It easily follows from these equations that $(s_{1},r_{1},\beta_{1})=(s_{2},r_{2},\beta_{2})=(1,0,0)$.

{\bf(1-ii)} If $b_{2}\not=0$ and $c_{2}=-\epsilon_{1}$, then we have that $b_{2}=2a_{2}+\epsilon_{1}$ by \eqref{eq2-4}.
Since $(a_{1},b_{1},c_{1})=(\epsilon_{1},0,2\epsilon_{1})$,
we have
\begin{eqnarray*}
\det A_{f}=(2\epsilon_{1}a_{2}+1)(b_{3}+c_{3}-2a_{3})=\pm 1
\end{eqnarray*}
Therefore, either {\bf(1-ii-a)}  $(a_{2},b_{2},c_{2})=(0,\epsilon_{1},-\epsilon_{1})$,  or  {\bf(1-ii-b)} $ (-\epsilon_{1},-\epsilon_{1},-\epsilon_{1})$ and
$b_{3}+c_{3}-2a_{3}=\pm 1$.

 {\bf(1-ii-a)} Suppose $(a_{2},b_{2},c_{2})=(0,\epsilon_{1},-\epsilon_{1})$, then
$a_{2}b_{3}-b_{2}a_{3}=-\epsilon_{1}a_{3}$.
As before,
if $a_{3}\not=0$ then $(s_{1},r_{1},\beta_{1})=(s_{2},r_{2},\beta_{2})=(1,0,0)$.
This satisfies the statement of proposition.
If $a_{3}=0$, then $b_{3}+c_{3}=\pm 1$ by the equation above.
From the relation $Z^{2}=-s_{1}XZ-r_{1}YZ-\beta_{1}XY$, we have
\begin{eqnarray}
& &-b_{3}^2=-s_{1}\epsilon_{1}b_{3}+r_{1}\epsilon_{1}b_{3}-\beta_{1}, \label{eqn;b3} \\
& &-c_{3}^2=s_{1}\epsilon_{1}c_{3}-r_{1}\epsilon_{1}c_{3}-\beta_{1}, \label{eqn;c3}\\
& &2b_{3}c_{3}=-2s_{1}\epsilon_{1}b_{3}-r_{1}\epsilon_{1}c_{3}+r_{1}\epsilon_{1}b_{3}-2\beta_{1}. \label{eqn;b3c3} 
\end{eqnarray}
 From these equations, we get
\begin{eqnarray*}
(b_{3}+c_{3})^2=1=-s_{1}\epsilon_{1}(b_{3}+c_{3}).
\end{eqnarray*}
Hence, $s_{1}=1$ and $b_{3}+c_{3}=-\epsilon_{1}$.
By \eqref{eqn;b3}, 
we have
\begin{eqnarray*}
-1+2\epsilon_{1}c_{3}-c_{3}^2=-\epsilon_{1}(-\epsilon_{1}-c_{3})+r_{1}\epsilon_{1}(-\epsilon_{1}-c_{3})-\beta_{1}.
\end{eqnarray*}
Substituting \eqref{eqn;c3} into this equation, we have
\begin{eqnarray*}
-1+2\epsilon_{1}c_{3}+\epsilon_{1}c_{3}-r_{1}\epsilon_{1}c_{3}-\beta_{1}=-\epsilon_{1}(-\epsilon_{1}-c_{3})+r_{1}\epsilon_{1}(-\epsilon_{1}-c_{3})-\beta_{1}.
\end{eqnarray*}
Hence,
\begin{eqnarray*}
2(2\epsilon_{1}c_{3}-1)=r_{1}=0.
\end{eqnarray*}
But this is impossible.
Therefore the case (1-ii-a) can not occur.

 {\bf(1-ii-b)} Suppose $(a_{2},b_{2},c_{2})=(-\epsilon_{1},-\epsilon_{1},-\epsilon_{1})$, then
$a_{2}b_{3}-b_{2}a_{3}=-\epsilon_{1}(b_{3}-a_{3})$.
With the method similar to that demonstrated  above, if $a_{3}\not=b_{3}$ then
$(s_{1},r_{1},\beta_{1})=(s_{2},r_{2},\beta_{2})=(1,0,0)$.
Hence, we may assume $a_{3}=b_{3}$.
Because $\det A_{f}=c_{3}+b_{3}-2a_{3}=\pm 1$, we also have $c_{3}-b_{3}=\pm 1$.
From the relation $Z^{2}=-s_{1}XZ-r_{1}YZ-\beta_{1}XY$, we have 
\begin{eqnarray}
& &b_{3}^{2}=-s_{1}\epsilon_{1}b_{3}+r_{1}\epsilon_{1}b_{3}+\beta_{1}, \label{eqn;2b3}\\
& &2b_{3}c_{3}-c_{3}^2=-2s_{1}\epsilon_{1}b_{3}+s_{1}\epsilon_{1}c_{3}+r_{1}\epsilon_{1}b_{3}+2\beta_{1}, \label{eqn;2b3c3-c3}\\
& &2b_{3}c_{3}=-2s_{1}\epsilon_{1}b_{3}+r_{1}\epsilon_{1}c_{3}+r_{1}\epsilon_{1}b_{3}+2\beta_{1}.\label{eqn;2b3c3}
\end{eqnarray}
By using \eqref{eqn;2b3c3-c3}  and \eqref{eqn;2b3c3}, we have 
\begin{eqnarray*}
c_{3}(r_{1}\epsilon_{1}-c_{3}-s_{1}\epsilon_{1})=0.
\end{eqnarray*}
Therefore, we have either  
$c_{3}=0$, or $c_{3}\not=0$ and $r_{1}\epsilon_{1}-c_{3}-s_{1}\epsilon_{1}=0$, i.e., $c_{3}=\epsilon_{1}(r_{1}-s_{1})$ with $r_{1}\not=s_{1}$.

We claim $c_{3}\not=0$.
If $c_{3}=0$, then
by using $\det A_{f}=\pm 1$ and $a_{3}=b_{3}$, we may put $b_{3}=\epsilon_{2}$.
By using \eqref{eqn;2b3c3-c3}  and \eqref{eqn;2b3c3} again, we have that
\begin{eqnarray*}
-2s_{1}\epsilon_{1}\epsilon_{2}+r_{1}\epsilon_{1}\epsilon_{2}+2\beta_{1}=0.
\end{eqnarray*}
Hence, it is easy to check that $(s_{1},r_{1},\beta_{1})=(0,0,0)$ or $(1,0,\epsilon_{1}\epsilon_{2})$.
However, using \eqref{eqn;2b3}, both of the cases give contradictions.
Consequently, $c_{3}\not=0$, i.e., $c_{3}=\epsilon_{1}(r_{1}-s_{1})$ with $r_{1}\not=s_{1}$.

Because $r_{1}\not=s_{1}$, there are two cases: $(s_{1},r_{1})=(1,0)$ and $(0,1)$.
We first assume that $(s_{1},r_{1})=(1,0)$.
In this case, $c_{3}=-\epsilon_{1}$.
By using \eqref{eqn;2b3c3-c3}, we have  $\beta_{1}=0$.
Therefore, this case gives $(s_{1},r_{1},\beta_{1})=(s_{2},r_{2},\beta_{2})=(1,0,0)$.
We next assume that $(s_{1},r_{1})=(0,1)$.
In this case, $c_{3}=\epsilon_{1}$.
Similarly, we have that $\epsilon_{1}b_{3}-1=2\beta_{1}$.
This also gives the equation
\begin{eqnarray*}
\epsilon_{1}b_{3}-1=\epsilon_{1}(b_{3}-\epsilon_{1})=2\beta_{1}.
\end{eqnarray*}
Recall that $b_{3}-c_{3}=\pm 1$ and $c_{3}=\epsilon_{1}$.
This gives a contradiction.
This finishes Case 1.


\textbf{Case 2: $c_{1}=0$.}
In this case we divided into two sub-cases: {\bf(2-i)} $(a_{1},b_{1},c_{1})=(\epsilon_{1},0,0)$,  and  {\bf(2-ii)} 
$(a_{1},b_{1},c_{1})=(\epsilon_{1},2\epsilon_{1},0)$.

 {\bf(2-i)} Assume $(a_{1},b_{1},c_{1})=(\epsilon_{1},0,0)$.
Then, it follows from \eqref{eq2-4}, \eqref{eq2-5} and \eqref{eq2-6} that
\begin{eqnarray}
& &\label{eq2-4'}2a_{2}b_{2}-b_{2}^2-c_{2}^2\beta_{2}=-\epsilon_{1}b_{2}; \\
& &\label{eq2-5'}2a_{2}c_{2}-c_{2}^2s_{2}=-\epsilon_{1}c_{2}; \\
& &\label{eq2-6'}2b_{2}c_{2}-c_{2}^2r_{2}=0.
\end{eqnarray}
By \eqref{eq2-5'} and \eqref{eq2-6'}, either 
{\bf(2-i-a)}  $c_{2}\not=0$ and $2a_{2}=c_{2}s_{2}-\epsilon_{1}$, $2b_{2}=c_{2}r_{2}$, or
 {\bf(2-i-b)}  $c_{2}=0$.

{\bf(2-i-a)}  First assume $c_{2}\not=0$.
Then, by $2a_{2}=c_{2}s_{2}-\epsilon_{1}$, we have $s_{2}=1$ and $c_{2}=2a_{2}+\epsilon_{1}$.
By substituting this equation into \eqref{eq2-6'},
we have that $r_{2}=0=b_{2}$.
Hence, by \eqref{eq2-4'}, $\beta_{2}=0$, i.e., $(s_{2},r_{2},\beta_{2})=(1,0,0)$.
Because $\det A_{f}=\pm 1$, we may put $b_{3}=\epsilon_{2}$.
Moreover, we have $\det A_{f}=-\epsilon_{1}\epsilon_{2}(2a_{2}+\epsilon_{1})=\pm 1$; therefore,
$a_{2}=0$ or $-\epsilon_{1}$.
If $a_{2}=-\epsilon_{1}$, then $a_{2}b_{3}-a_{3}b_{2}=-\epsilon_{1}\epsilon_{2}\not=0$.
Hence, with the method similar to that demonstrated in Case 1, we have $(s_{1},r_{1},\beta_{1})=(s_{2},r_{2},\beta_{2})=(1,0,0)$.
Thus, we may assume $a_{2}=0$, i.e.,
\begin{eqnarray*}
A_{f}=
\left(
\begin{array}{ccc}
\epsilon_{1} & 0 & 0 \\
0 & 0 & \epsilon_{1} \\
a_{3} & \epsilon_{2} & c_{3}
\end{array}
\right).
\end{eqnarray*}
By using $Z^{2}=-s_{1}XZ-r_{1}YZ-\beta_{1}XY$ and $(s_{2},r_{2},\beta_{2})=(1,0,0)$,
it is easy to get that
\begin{eqnarray*}
& &2a_{3}\epsilon_{2}-1=-s_{1}\epsilon_{1}\epsilon_{2}; \\
& &2\epsilon_{2}c_{3}=-r_{1}\epsilon_{1}\epsilon_{2}; \\
& &(2a_{3}-c_{3})c_{3}=-s_{1}\epsilon_{1}c_{3}-r_{1}\epsilon_{1}a_{3}+r_{1}c_{3}\epsilon_{1}-\beta_{1}.
\end{eqnarray*}
By using the first and second equations, we have $s_{1}=1$, $r_{1}=0$ and $c_{3}=0$.
Therefore, by the third equation, we have that $(s_{1},r_{1},\beta_{1})=(s_{2},r_{2},\beta_{2})=(1,0,0)$.
Consequently, if $(a_{1},b_{1},c_{1})=(\epsilon_{1},0,0)$ and $c_{2}\not=0$,
then $(s_{1},r_{1},\beta_{1})=(s_{2},r_{2},\beta_{2})=(1,0,0)$.

{\bf(2-i-b)}  We next assume $c_{2}=0$.
Because $\det A_{f}=\epsilon_{1}b_{2}c_{3}=\pm 1$, we may put $b_{2}=\epsilon_{2}$ and $c_{3}=\epsilon_{3}$,
i.e.,
\begin{eqnarray*}
A_{f}=
\left(
\begin{array}{ccc}
\epsilon_{1} & 0 & 0 \\
a_{2} & \epsilon_{2} & 0 \\
a_{3} & b_{3} & \epsilon_{3}
\end{array}
\right).
\end{eqnarray*}
Then, it follows from \eqref{eq2-4'} that
$2a_{2}\epsilon_{2}-1=-\epsilon_{1}\epsilon_{2}$, i.e., $a_{2}=\frac{-\epsilon_{1}+\epsilon_{2}}{2}$.
By using $Z^{2}=-s_{1}XZ-r_{1}YZ-\beta_{1}XY$, 
it is easy to get that
\begin{eqnarray*}
& &2a_{3}b_{3}-b_{3}^2-\beta_{2}=-s_{1}\epsilon_{1}b_{3}-r_{1}(a_{2}b_{3}+a_{3}\epsilon_{2}-\epsilon_{2}b_{3})-\beta_{1}\epsilon_{1}\epsilon_{2}; \\
& &2b_{3}\epsilon_{3}-r_{2}=-r_{1}\epsilon_{2}\epsilon_{3}; \\
& &2a_{3}\epsilon_{3}-s_{2}=-s_{1}\epsilon_{1}\epsilon_{3}-r_{1}a_{2}\epsilon_{3}.
\end{eqnarray*}

If $\epsilon_{1}=\epsilon_{2}$, then $a_{2}=0$ and
\begin{eqnarray*}
& &2a_{3}b_{3}-b_{3}^2-\beta_{2}=-s_{1}\epsilon_{1}b_{3}-r_{1}(a_{3}\epsilon_{1}-\epsilon_{1}b_{3})-\beta_{1}; \\
& &2b_{3}\epsilon_{3}-r_{2}=-r_{1}\epsilon_{1}\epsilon_{3}; \\
& &2a_{3}\epsilon_{3}-s_{2}=-s_{1}\epsilon_{1}\epsilon_{3}.
\end{eqnarray*}
By using the second and third equations, we have that $(s_{1},r_{1})=(s_{2},r_{2})$.
Therefore, if $\epsilon_{1}=\epsilon_{3}$, then we also have $b_{3}=a_{3}=0$.
Using the first equation, we have $\beta_{1}=\beta_{2}$, i.e., $(s_{1},r_{1},\beta_{1})=(s_{2},r_{2},\beta_{2})$.
Suppose $\epsilon_{1}\not=\epsilon_{3}$, i.e., $\epsilon_{3}=-\epsilon_{1}$.
In this case, if $s_{1}=s_{2}=0$ (resp.\ $s_{1}=s_{2}=1$)
then $a_{3}=0$ (resp.\ $a_{3}=-\epsilon_{1}$) by using the third equation.
Similarly by using the second equation, if $r_{1}=r_{2}=0$ (resp.\ $r_{1}=r_{2}=1$)
then $b_{3}=0$ (resp.\ $b_{3}=-\epsilon_{1}$).
Therefore, by using the first equation, it is easy to check that
$\beta_{1}=\beta_{2}$.
Consequently, in the case when $\epsilon_{1}=\epsilon_{2}$, hence $(a_{2},b_{2},c_{2})=(0,\epsilon_{1},0)$, we have $(s_{1},r_{1},\beta_{1})=(s_{2},r_{2},\beta_{2})$, i.e.,
this case satisfies the statement of proposition.

If $-\epsilon_{1}=\epsilon_{2}$, then $a_{2}=-\epsilon_{1}$ and
\begin{eqnarray*}
& &2a_{3}b_{3}-b_{3}^2-\beta_{2}=-s_{1}\epsilon_{1}b_{3}+r_{1}a_{3}\epsilon_{1}+\beta_{1}; \\
& &2b_{3}\epsilon_{3}-r_{2}=r_{1}\epsilon_{1}\epsilon_{3}; \\
& &2a_{3}\epsilon_{3}-s_{2}=-s_{1}\epsilon_{1}\epsilon_{3}+r_{1}\epsilon_{1}\epsilon_{3}.
\end{eqnarray*}
By using the second equation, we have that $r_{1}=r_{2}$.
If $r_{1}=r_{2}=0$, then $b_{3}=0$ by the second equation and
$s_{1}=s_{2}$ by the third equation.
Moreover, by using the first equation, we have $(s_{1},0,\beta_{1})=(s_{2},0,-\beta_{2})$.
This implies that $(3)-(a)$ and $(3)-(b)$ in the statement of the proposition.
If $r_{1}=r_{2}=1$, then
$b_{3}=\frac{\epsilon_{1}+\epsilon_{3}}{2}$ by the second equation and
$s_{1}\not=s_{2}$ by the third equation.
We first assume $(s_{1},s_{2})=(1,0)$.
Then, by the third equation, we have that $a_{3}=0$.
Therefore, the first equation gives
\begin{eqnarray*}
-\frac{1+\epsilon_{1}\epsilon_{3}}{2}-\beta_{2}=-\frac{1+\epsilon_{1}\epsilon_{3}}{2}+\beta_{1}.
\end{eqnarray*}
Therefore, $\beta_{1}=-\beta_{2}$, i.e.,
$(s_{1},r_{1},\beta_{1})$ and $(s_{2},r_{2},\beta_{2})$ are the pair
$(1,1,r)$ and $(0,1,-r)$.
This implies that $(3)-(c)$ in the statement of the proposition.
We next assume $(s_{1},s_{2})=(0,1)$.
Then, by the second and third equations, we have that $a_{3}=b_{3}$.
Therefore, the first equation gives
\begin{eqnarray*}
\frac{1+\epsilon_{1}\epsilon_{3}}{2}-\beta_{2}=\frac{1+\epsilon_{1}\epsilon_{3}}{2}+\beta_{1}.
\end{eqnarray*}
Therefore, $\beta_{1}=-\beta_{2}$, i.e.,
$(s_{1},r_{1},\beta_{1})$ and $(s_{2},r_{2},\beta_{2})$ are the pair
$(0,1,r)$ and $(1,1,-r)$.
This implies that $(3)-(c)$ in the statement of the proposition.
Consequently, if $(a_{1},b_{1},c_{1})=(\epsilon_{1},0,0)$ and $c_{2}=0$, then the statement holds.
Therefore the first sub-case  {\bf(2-i)} is done.

{\bf(2-ii)}  Assume $(a_{1},b_{1},c_{1})=(\epsilon_{1},2\epsilon_{1},0)$.
Then, it follows from \eqref{eq2-4}, \eqref{eq2-5} and \eqref{eq2-6} that
\begin{eqnarray}
& &\label{eq2-4''}2a_{2}b_{2}-b_{2}^2-c_{2}^2\beta_{2}=\epsilon_{1}b_{2}-2\epsilon_{1}a_{2}; \\
& &\label{eq2-5''}2a_{2}c_{2}-c_{2}^2s_{2}=-\epsilon_{1}c_{2}; \\
& &\label{eq2-6''}2b_{2}c_{2}-c_{2}^2r_{2}=-2\epsilon_{1}c_{2}.
\end{eqnarray}
By \eqref{eq2-5''} and \eqref{eq2-6''}, either 
{\bf(2-ii-a)}  $c_{2}\not=0$ and $2a_{2}=c_{2}s_{2}-\epsilon_{1}$, $2b_{2}=c_{2}r_{2}-2\epsilon_{1}$, or
{\bf(2-ii-b)}  $c_{2}=0$.

{\bf(2-ii-a)}  We first assume $c_{2}\not=0$.
Then, by $2a_{2}=c_{2}s_{2}-\epsilon_{1}$, we have  $s_{2}=1$ and $c_{2}=2a_{2}+\epsilon_{1}$.
Substituting this equation into $2b_{2}=c_{2}r_{2}-2\epsilon_{1}$, we have  $r_{2}=0$ and $b_{2}=-\epsilon_{1}$.
Therefore, $\beta_{2}=0$ by \eqref{eq2-4''}.
By using $Z^{2}=-s_{1}XZ-r_{1}YZ-\beta_{1}XY$ and $(s_{2},r_{2},\beta_{2})=(1,0,0)$, 
it is easy to get that 
\begin{eqnarray}
& &2a_{3}b_{3}-b_{3}^2=-s_{1}(-\epsilon_{1}b_{3}+2\epsilon_{1}a_{3})
-r_{1}(a_{2}b_{3}-\epsilon_{1}a_{3}+\epsilon_{1}b_{3})
-\beta_{1}(1+2\epsilon_{1}a_{2}); \label{eqn;a3b3c3-1}\\
& &2b_{3}c_{3}= -2s_{1}\epsilon_{1}c_{3} -r_{1}(-\epsilon_{1}c_{3}+2a_{2}b_{3}+\epsilon_{1}b_{3})-\beta_{1}(4a_{2}\epsilon_{1}+2); \label{eqn;a3b3c3-2}\\
& &(2a_{3}-c_{3})c_{3}=-s_{1}\epsilon_{1}c_{3}-r_{1}(-a_{2}c_{3}+2a_{2}a_{3}+\epsilon_{1}a_{3}-\epsilon_{1}c_{3})-\beta_{1}(2a_{2}\epsilon_{1}+1).\label{eqn;a3b3c3-3}
\end{eqnarray}

Because $\det A_{f}=(2a_{2}\epsilon_{1}+1)(2a_{3}-b_{3}-c_{3})=\pm 1$, either {\bf(2-ii-a-I)}  $a_{2}=0$ or {\bf(2-ii-a-II)}  $a_{2}=-\epsilon_{1}$,
and we may put
$2a_{3}-b_{3}-c_{3}=\epsilon_{3}$.

{\bf(2-ii-a-I)}  Assume $a_{2}=0$.
With the method similar to that demonstrated in Case 1, if $a_{2}b_{3}-a_{3}b_{2}=a_{3}\not=0$ then $(s_{1},r_{1},\beta_{1})=(s_{2},r_{2},\beta_{2})=(1,0,0)$.
Therefore, we may assume $a_{3}=0$ and $-b_{3}-c_{3}=\epsilon_{3}$.
Hence, by the above equations, we have that
\begin{eqnarray}
& &-b_{3}^2=s_{1}\epsilon_{1}b_{3}-r_{1}\epsilon_{1}b_{3}-\beta_{1}; \label{eqn;a3b3c3'-1}\\
& &2b_{3}c_{3}= -2s_{1}\epsilon_{1}c_{3} -r_{1}(-\epsilon_{1}c_{3}+\epsilon_{1}b_{3})-2\beta_{1}; \label{eqn;a3b3c3'-2}\\
& &-c_{3}^{2}=-s_{1}\epsilon_{1}c_{3}+r_{1}\epsilon_{1}c_{3}-\beta_{1}.\label{eqn;a3b3c3'-3}
\end{eqnarray}
This implies that
\begin{eqnarray*}
-(b_{3}+c_{3})^2=-1=s_{1}\epsilon_{1}(b_{3}+c_{3})=-s_{1}\epsilon_{1}\epsilon_{3}.
\end{eqnarray*}
Therefore, we have $s_{1}=1=\epsilon_{1}\epsilon_{3}$ and $c_{3}=-b_{3}-\epsilon_{1}$.
By substituting these equations into the third equation,
we have 
\begin{eqnarray*}
-b_{3}^{2}-2\epsilon_{1}b_{3}-1=\epsilon_{1}(b_{3}+\epsilon_{1})-r_{1}\epsilon_{1}(b_{3}+\epsilon_{1})-\beta_{1}.
\end{eqnarray*}
Because of the first equation,
we have 
\begin{eqnarray*}
2\epsilon_{1}b_{3}+2=r_{1}.
\end{eqnarray*}
This implies that $r_1=0$ and $b_3=-\epsilon_1$. Hence $c_3=-b_3-\epsilon_1=0$. Therefore,  from \eqref{eqn;a3b3c3'-2}, we have 
$\beta_1=0$.  
Therefore, $(s_{1},r_{1},\beta_{1})=(s_{2},r_{2},\beta_{2})=(1,0,0)$.
This satisfies the statement of proposition, and the case {\bf(2-ii-a-I)} is done.

{\bf(2-ii-a-II)} Assume $a_{2}=-\epsilon_{1}$
With the method similar to that demonstrated in Case 1, if $a_{3}\not=b_{3}$ then $(s_{1},r_{1},\beta_{1})=(s_{2},r_{2},\beta_{2})=(1,0,0)$.
Therefore, we may assume $a_{3}=b_{3}$ and $a_{3}-c_{3}=\epsilon_{3}$.
By the above equations \eqref{eqn;a3b3c3-1}, \eqref{eqn;a3b3c3-2}, and \eqref{eqn;a3b3c3-3}, we have
\begin{eqnarray*}
& &a_{3}^2=-s_{1}\epsilon_{1}a_{3}+r_{1}\epsilon_{1}a_{3}+\beta_{1}; \\
& &2a_{3}c_{3}= -2s_{1}\epsilon_{1}c_{3} -r_{1}(-\epsilon_{1}c_{3}-\epsilon_{1}a_{3})+2\beta_{1}; \\
& &(2a_{3}-c_{3})c_{3}=-s_{1}\epsilon_{1}c_{3}+r_{1}\epsilon_{1}a_{3}+\beta_{1}.
\end{eqnarray*}
This implies that
\begin{eqnarray*}
(a_{3}+c_{3})(-a_{3}+c_{3})&=&s_{1}\epsilon_{1}a_{3}-s_{1}\epsilon_{1}c_{3} +r_{1}\epsilon_{1}c_{3}-r_{1}\epsilon_{1}a_{3} \\
&=&\epsilon_{1}(r_{1}-s_{1})(-a_{3}+c_{3}).
\end{eqnarray*}
Because $a_{3}-c_{3}=\epsilon_{3}$, we have that $a_{3}+c_{3}=\epsilon_{1}(r_{1}-s_{1})$; therefore, $r_{1}\not=s_{1}$.
If $(s_{1},r_{1})=(0,1)$, then $2a_{3}c_{3}=1+2\beta_{1}$
by the second equation above.
This gives a contradiction.
Hence, $(s_{1},r_{1})=(1,0)$.
In this case, $a_{3}=\frac{-\epsilon_{1}+\epsilon_{3}}{2}$ and $c_{3}=\frac{-\epsilon_{1}-\epsilon_{3}}{2}$.
If $\epsilon_{1}=\epsilon_{3}$, then $a_{3}=0$ and $c_{3}=-\epsilon_{1}$.
In this case, by using the first equation, $\beta_{1}=0$.
However, by using the second equation, we also have $\beta_{1}=-1$.
This gives a contradiction and we have $\epsilon_{1}=-\epsilon_{3}$, i.e., $a_{3}=-\epsilon_{1}$ and $c_{3}=0$.
It is easy to check that $(s_{1},r_{1},\beta_{1})=(s_{2},r_{2},\beta_{2})=(1,0,0)$.
Consequently, if $(a_{1},b_{1},c_{1})=(\epsilon_{1},2\epsilon_{1},0)$ and $c_{2}\not=0$,
then $(s_{1},r_{1},\beta_{1})=(s_{2},r_{2},\beta_{2})=(1,0,0)$.
This satisfies the statement of proposition.
This finishes the proof for {\bf(2-ii-a)}.

{\bf(2-ii-b)} We next assume $c_{2}=0$, i.e.,
\begin{eqnarray*}
A_{f}=
\left(
\begin{array}{ccc}
\epsilon_{1} & 2\epsilon_{1} & 0 \\
a_{2} & b_{2} & 0 \\
a_{3} & b_{3} & c_{3}
\end{array}
\right).
\end{eqnarray*}
Since $\det  A_f=\pm1$, we have $c_3=\pm1=:\epsilon_3$.
By \eqref{eq2-4''}
\begin{eqnarray*}
2a_{2}b_{2}-b_{2}^2=\epsilon_{1}b_{2}-2\epsilon_{1}a_{2}.
\end{eqnarray*}
Hence,
\begin{eqnarray*}
(2a_{2}-b_{2})(b_{2}+\epsilon_{1})=0.
\end{eqnarray*}
Therefore,
$b_{2}=2a_{2}$ or $-\epsilon_{1}$.
If $b_{2}=2a_{2}$,
then $\det A_f=0$, which is  contradiction.
Therefore, $b_{2}=-\epsilon_{1}$.
Hence, $\det A_{f}=\epsilon_{3}(-1-2\epsilon_{1}a_{2})=\pm 1$; therefore,
\begin{itemize}
\item $a_{2}=0$ or
\item $a_{2}=-\epsilon_{1}$.
\end{itemize}
By using $Z^{2}=-s_{1}XZ-r_{1}YZ-\beta_{1}XY$, 
it is easy to get that 
\begin{eqnarray*}
& &2a_{3}b_{3}-b_{3}^2-\beta_{2}=-s_{1}(-\epsilon_{1}b_{3}+2\epsilon_{1}a_{3})
-r_{1}(a_{2}b_{3}-\epsilon_{1}a_{3}+\epsilon_{1}b_{3})
-\beta_{1}(1+2\epsilon_{1}a_{2}); \\
& &2b_{3}\epsilon_{3}-r_{2}
=-2s_{1}\epsilon_{1}\epsilon_{3}+r_{1}\epsilon_{1}\epsilon_{3}; \\
& &2a_{3}-s_{2}\epsilon_{3}=-s_{1}\epsilon_{1}-r_{1}a_{2}.
\end{eqnarray*}
By the second equation, we have that $r_{1}=r_{2}$.
If $r_{1}=r_{2}=0$, by the second and third equations, we have that $b_{3}=-s_{1}\epsilon_{1}$ and $s_{1}=s_{2}$, respectively.
It follows easily from the first equation that $\beta_{1}=\beta_{2}$ for $a_{2}=0$ and $\beta_{1}=-\beta_{2}$ for $a_{2}=-\epsilon_{1}$.
This implies that $(3)-(a)$ and $(3)-(b)$ and $(3)$ with $(s_{1},0,\beta_{1})=(s_{2},0,\beta_{2})$ in the statement of the proposition.
If $r_{1}=r_{2}=1$, then by the above equations, we have that
\begin{eqnarray*}
& &2a_{3}b_{3}-b_{3}^2-\beta_{2}=-s_{1}(-\epsilon_{1}b_{3}+2\epsilon_{1}a_{3})
-a_{2}b_{3}+\epsilon_{1}a_{3}-\epsilon_{1}b_{3}
-\beta_{1}(1+2\epsilon_{1}a_{2}); \\
& &2b_{3}\epsilon_{3}-1
=-2s_{1}\epsilon_{1}\epsilon_{3}+\epsilon_{1}\epsilon_{3}; \\
& &2a_{3}-s_{2}\epsilon_{3}=-s_{1}\epsilon_{1}-a_{2}.
\end{eqnarray*}

When $a_{2}=0$, then by the third equation we have that $s_{1}=s_{2}$.
If $s_{1}=s_{2}=0$, then by the third equation we have $a_{3}=0$; therefore by the first and second equations we have
\begin{eqnarray*}
-\frac{1+\epsilon_{1}\epsilon_{3}}{2}-\beta_{2}=-\frac{1+\epsilon_{1}\epsilon_{3}}{2}-\beta_{1}
\end{eqnarray*}
Hence, $\beta_{1}=\beta_{2}$.
This implies that $(3)$ with $(0,1,\beta_{1})=(0,1,\beta_{2})$ in the statement of the proposition.
If $s_{1}=s_{2}=1$, then by the second and third equations, we have that $a_{3}=b_{3}=\frac{-\epsilon_{1}+\epsilon_{3}}{2}$.
Using the first equation, we have
$\beta_{1}=\beta_{2}$.
This implies that $(3)$ with $(1,1,\beta_{1})=(1,1,\beta_{2})$ in the statement of the proposition.

When $a_{2}=-\epsilon_{1}$, then by the third equation we have that $s_{1}\not=s_{2}$.
If $(s_{1},s_{2})=(1,0)$, then it follows from the third equation that $a_{3}=0$; therefore by the first and second equations we have
\begin{eqnarray*}
-\frac{1-\epsilon_{1}\epsilon_{3}}{2}-\beta_{2}=-\frac{1-\epsilon_{1}\epsilon_{3}}{2}+\beta_{1}
\end{eqnarray*}
Hence, $\beta_{1}=-\beta_{2}$.
If $(s_{1},s_{2})=(0,1)$, then
by the second and third equations, we have that $a_{3}=b_{3}=\frac{\epsilon_{1}+\epsilon_{3}}{2}$.
Using the first equation, we have
$\beta_{1}=-\beta_{2}$.
This implies that $(3)-(c)$ in the statement of the proposition.
Consequently, if $(a_{1},b_{1},c_{1})=(\epsilon_{1},2\epsilon_{1},0)$ and $c_{2}=0$, then the statement holds.
Therefore {\bf(2-ii-b)} is finished, and 
this establishes the statement $(2)\Rightarrow (3)$.

Finally, we prove $(3)\Rightarrow (1)$.
If $(s_{1},r_{1},\beta_{1})=(s_{2},r_{2},\beta_{2})$, then the statement is trivial.
Assume $(s_{1},r_{1},\beta_{1})\not=(s_{2},r_{2},\beta_{2})$.
Recall that $H_{1}\cong S^3\times_{T^1}P(\C_{1}\oplus\C)$. Let $f:H_{1}\to H_{1}$ be the diffeomorphism
which is induced from the composition of the diffeomorphisms
\begin{eqnarray*}
S^{3}\times_{T^1}P(\C_{1}\oplus \C)\stackrel{g}{\to} S^{3}\times_{T^1}P(\C_{-1}\oplus \C)\stackrel{h}{\to} S^{3}\times_{T^1}P(\C_{1}\oplus \C),
\end{eqnarray*}
where $g$ is the diffeomorphism induced from the orientation reversing of the fibers and
$h$ is the diffeomorphism induced from the tensor product of the tautological line bundle on $\gamma_{-1}\oplus\epsilon$.
Then, it is easy to check that the induced
homomorphism $f^{*}$ is $f^{*}(X)=x$ and $f^{*}(Y)=-x-y$,
where $H^{*}(H_{1})\simeq \Z[x,y]/\langle x^2, y^2+xy \rangle$.
Then, we can easily check the following isomorphisms;
\begin{eqnarray*}
& &f^{*}\xi_{(0,0,\beta)}\equiv \xi_{(0,0,-\beta)}; \\
& &f^{*}\xi_{(1,0,\beta)}\equiv \xi_{(1,0,-\beta)}; \\
& &f^{*}\xi_{(0,1,\beta)}\equiv \xi_{(-1,-1,-\beta)}.
\end{eqnarray*}
Because of Lemma \ref{formula of Chern class}, we have 
\begin{eqnarray*}
\gamma_{x+y}\otimes \xi_{(-1,-1,-\beta)}\equiv \xi_{(1,1,-\beta)},
\end{eqnarray*}
where $\gamma_{x+y}$ is the line bundle over $H_{1}$ induced from $x+y\in H^{2}(H_{1})$.
This establishes that
\begin{eqnarray*}
& &P(\xi_{(0,0,\beta)}) \cong P(\xi_{(0,0,-\beta)}); \\
& &P(\xi_{(1,0,\beta)}) \cong P(\xi_{(1,0,-\beta)}); \\
& &P(\xi_{(0,1,\beta)}) \cong P(\xi_{(1,1,-\beta)}).
\end{eqnarray*}
\end{proof}

Consequently, using Theorem  \ref{2-stage} and \ref{3-stage},
we have Theorem \ref{main}.

\section{Cohomological non-rigidity of $8$-dimensional $\C P$-tower}
\label{sect5}

In this section, we classify all $2$-stage $\C P$-towers whose first stage is $\C P^{3}$.
We first introduce the following classification result of complex $2$-dimensional  vector bundles over $\C P^3$  by Atiyah and Rees \cite{AtRe}.
Let ${\rm Vect}_2(\C P^3)$ be the set of complex  $2$-dimensional vector bundles over $\C P^3$ up to bundle isomorphisms.
\begin{theorem}[Atiyah-Rees]
\label{Atiyah-Rees}
There exist an injective map $\phi: {\rm Vect}_2(\C P^3)\to \Z_2\oplus\Z\oplus\Z$ such that $\phi(\xi)=(\alpha(\xi), c_1(\xi), c_2(\xi))$,
where $c_1(\xi)$ and $c_2(\xi)$ are the first and the second Chern classes of $\xi$,  and $\alpha(\xi)$ is a mod $2$ element which is $0$ 
when $c_1(\xi)$ is odd.
\end{theorem}
By Theorem~\ref{Atiyah-Rees}, any element in ${\rm Vect}_2(\C P^3)$ can be denoted by $\eta_{(\alpha,c_{1},c_{2})}$, where
$(\alpha,c_{1},c_{2})\in \Z_{2}\oplus \Z\oplus \Z$ such that $\alpha\equiv 0$ (mod $2$) when $c_{1}\equiv1$ (mod $2$).
The goal of this section is to classify the topological types of $P(\eta_{(\alpha,c_{1},c_{2})})$ up to diffeomorphisms.

Because $P(\eta_{(\alpha,c_{1},c_{2})})$ is diffeomorphic to $ P(\eta_{(\alpha,c_{1},c_{2})}\otimes \gamma)$ for any line bundle $\gamma$ over
$\C P^3$ by Lemma \ref{tensor line budle},
we may assume $c_{1}\in \{0,1\}$.
Therefore, in order to classify all  $P(\eta_{(\alpha,c_{1},c_{2})})$ up to diffeomorphisms, it is enough to classify
the following:
\begin{eqnarray*}
M_{0}(u)=P(\eta_{(0,0,u)}); \\
M_{1}(u)=P(\eta_{(1,0,u)}); \\
N(u)=P(\eta_{(0,1,u)}),
\end{eqnarray*}
where $u\in \Z$.
In the following three lemmas, we classify the cohomology rings of the above three types of manifolds up to graded ring isomorphisms.
\begin{lemma}
\label{cohom-8dim-1}
Two cohomology rings $H^{*}(M_{\alpha}(u))$ and $H^{*}(N(u'))$ are not isomorphic for any $u,\ u'\in \Z$.
\end{lemma}
\begin{proof}
By the Borel-Hirzebruch formula \eqref{BH-formula}, we have ring isomorphisms
\begin{eqnarray*}
& &H^{*}(M_{\alpha}(u))\simeq \Z[X,Y]/\langle X^{4},\ uX^{2}+Y^{2} \rangle, \  \textrm{ and} \\
& &H^{*}(N(u'))\simeq \Z[x,y]/\langle x^{4},\ u'x^{2}+xy+y^{2} \rangle.
\end{eqnarray*}
Assume that there is an isomorphism map $f:H^{*}(M_{\alpha}(u))\to H^{*}(N(u'))$.
Then  we may put
\begin{eqnarray*}
& &f(X)=ax+by, \  \textrm{ and} \\
& &f(Y)=cx+dy,
\end{eqnarray*}
for some $a,\ b,\ c,\ d\in \Z$ such that $ad-bc=\epsilon=\pm 1$.
By taking the inverse of $f$, we also have 
\begin{eqnarray*}
& &f^{-1}(x)=d\epsilon X-b\epsilon Y, \  \textrm{ and}\\
& &f^{-1}(y)=-c\epsilon X+a\epsilon Y.
\end{eqnarray*}

From  the ring structures of $H^{*}(M_{\alpha}(u))$ and $H^{*}(N(u'))$, we have
$f(uX^{2}+Y^{2})=0$ and $f^{-1}(y^{2}+xy+u'x^{2})=0$.
Therefore 
we have the following equations:
\begin{eqnarray}
& &\label{1st-8,0} u(a^{2}-u'b^{2})+(c^{2}-u'd^{2})=0; \\
& &\label{2nd-8,0} u(2ab-b^{2})+(2cd-d^{2})=0; \\
& &\label{3rd-8,0} c^{2}-a^{2}u -cd+abu +u'd^{2}-b^{2}uu'=0; \\
& &\label{4th-8,0} -2ac +cb+ad -2bdu'=0.
\end{eqnarray}
Because $f^{-1}(x^{4})=(d X-b Y)^{4}=0$, we also have
\begin{eqnarray*}
bd(d^{2}-ub^{2})=0.
\end{eqnarray*}
Therefore $bd=0$, or otherwise $d^{2}=ub^{2}$.
We first assume $bd=0$.
Then, there are two cases: $b=0$ and $d=0$.
If $b=0$, then $|a|=|d|=1$.
However, by using \eqref{2nd-8,0}, we have $2cd=1$.
This gives a contradiction.
If $d=0$, then $|b|=|c|=1$.
By using \eqref{4th-8,0}, we have $c(-2a+b)=0$, i.e., $b=2a$ by $|c|=1$.
However, this contradicts to $|b|=1$.
Hence, $bd\not=0$ and $d^{2}=ub^{2}$, i.e., $|d|=\sqrt{|u|} |b|$.
In this case, because $ad-bc=\epsilon=\pm 1$, we have $|b|=1$ and $d^{2}=u$.
Let $b=\epsilon'=\pm 1$ and $d=\sqrt{u} \epsilon''$, where $\epsilon''=\pm 1$.
Then, it follows from $ad-bc=\epsilon$ that $c=-\epsilon\epsilon'+a\sqrt{u} \epsilon''\epsilon'$.
Therefore, by using \eqref{1st-8,0}, we have the following equation:
\begin{eqnarray*}
& &u(a^{2}-u'b^{2})+(c^{2}-u'd^{2}) \\
&=&u(a^{2}-u')+(-\epsilon\epsilon'+a\sqrt{u} \epsilon''\epsilon')^{2}-u'u \\
&=&2ua^{2}-2uu'+1-2a\sqrt{u}\epsilon \epsilon''=0.
\end{eqnarray*}
However, this gives the equation $1=2(-ua^{2}+uu'+a\sqrt{u}\epsilon \epsilon'')$, which is a contradiction.
Hence, $H^{*}(M_{\alpha}(u))\not\simeq H^{*}(N(u'))$ for all $u,\ u'\in \Z$.
\end{proof}

\begin{lemma}
\label{cohom-8dim-2}
The following two statements are equivalent.
\begin{enumerate}
\item $H^{*}(M_{\alpha}(u))\simeq H^{*}(M_{\alpha'}(u'))$ where $\alpha,\ \alpha'\in \{0,1\}$.
\item $u=u'\in \Z$
\end{enumerate}
\end{lemma}
\begin{proof}
Because $(2)\Rightarrow (1)$ is trivial, it is enough to show $(1)\Rightarrow (2)$. Assume there is an
isomorphism
$f: H^{*}(M_{\alpha}(u))\simeq H^{*}(M_{\alpha'}(u'))$  where
\begin{eqnarray*}
& &H^{*}(M_{\alpha}(u))\simeq \Z[X,Y]/\langle X^{4},\ uX^{2}+Y^{2} \rangle; \\
& &H^{*}(M_{\alpha'}(u'))\simeq \Z[x,y]/\langle x^{4},\ u'x^{2}+y^{2} \rangle.
\end{eqnarray*}
We may use the same representation for $f$ as in the proof of Lemma~\ref{cohom-8dim-1}.
Note that
 $f(uX^{2}+Y^{2})=0$ and $f^{-1}(u'x^{2}+y^{2})=0$.
By using the representation of $f$,
we have the following equations:
\begin{eqnarray}
& &\label{1st-8,1} ua^{2}-uu'b^{2}+c^{2}-u'd^{2}=0; \\
& &\label{2nd-8,1} uab+cd=0; \\
& &\label{3rd-8,1} u'd^{2}-uu'b^{2}+c^{2}-a^{2}u=0; \\
& &\label{4th-8,1} u'bd+ac=0.
\end{eqnarray}
By \eqref{1st-8,1} and \eqref{3rd-8,1}, we have
\begin{eqnarray}
& &\label{5th-8,1}c^{2}=b^{2}uu'; \\
& &\label{6th-8,1}ua^{2}=u'd^{2}.
\end{eqnarray}
Because $X^{4}=0$, we also have that
\begin{eqnarray*}
& &ab(a^{2}-b^{2}u')=0.
\end{eqnarray*}
We first assume $ab\not=0$.
Then
\begin{eqnarray*}
a^{2}=b^{2}u'
\end{eqnarray*}
 by this equation.
Together with \eqref{5th-8,1} and \eqref{6th-8,1}, we have that
\begin{eqnarray*}
c^{2}b^{2}=b^{4}uu'=b^{2}a^{2}u=b^{2}d^{2}u'=a^{2}d^{2}.
\end{eqnarray*}
This implies that
\begin{eqnarray*}
(ad-bc)(ad+bc)=\epsilon(ad+bc)=0.
\end{eqnarray*}
Hence, $ad=-bc$. However this gives a contradiction because $ad-bc=2ad=\epsilon=\pm 1$.
Consequently, we have $ab=0$.
Since $ad-bc=\epsilon$, if $a=0$ then $|b|=|c|=1$; therefore, we have $u=u'=\pm 1$ by \eqref{5th-8,1};
if $b=0$ then $|a|=|d|=1$; therefore, we have $u=u'$ by \eqref{6th-8,1}.
This establishes the statement.
\end{proof}

\begin{lemma}
\label{cohom-8dim-3}
The following two statements are equivalent.
\begin{enumerate}
\item $H^{*}(N(u))\simeq H^{*}(N(u'))$
\item $u=u'\in \Z$
\end{enumerate}
\end{lemma}
\begin{proof}
Because $(2)\Rightarrow (1)$ is trivial, it is enough to show $(1)\Rightarrow (2)$.
Assume there is an isomorphism $f: H^{*}(N(u))\simeq H^{*}(N(u'))$ where
\begin{eqnarray*}
& &H^{*}(N(u))\simeq \Z[X,Y]/\langle X^{4},\ uX^{2}+xy+Y^{2} \rangle; \\
& &H^{*}(N(u'))\simeq \Z[x,y]/\langle x^{4},\ u'x^{2}+xy+y^{2} \rangle.
\end{eqnarray*}

Again, we use the same representation for $f$ as in the proof of Lemma~\ref{cohom-8dim-1}.
Because $f(Y^{2}+XY+uX^{2})=0$ and $f^{-1}(y^{2}+xy+u'x^{2})=0$,
we have that
\begin{eqnarray}
& &\label{1st-8,2} c^{2}-d^{2}u'=-ua^{2}+b^{2}uu'-ac+bdu'; \\
& &\label{2nd-8,2} 2cd-d^{2}=-2abu+b^{2}u-ad-bc+bd; \\
& &\label{3rd-8,2} c^{2}-a^{2}u=-u'd^{2}+b^{2}uu'+cd-bau; \\
& &\label{4th-8,2} -2ac-a^{2}=2bdu'+b^{2}u'-ad-bc-ab.
\end{eqnarray}
Because $f(X^{4})=0$ and $f^{-1}(x^{4})=0$, there are the following two cases:
\begin{enumerate}
\item $b=0$;
\item $b\not=0$ and $4a^{3}-6a^{2}b+4ab^{2}(1-u')+b^{3}(2u'-1)=-4d^{3}-6d^{2}b-4db^{2}(1-u)+b^{3}(2u-1)=0$.
\end{enumerate}

If $b=0$, then $|a|=|d|=1$.
Therefore, by \eqref{2nd-8,2}, $2c=d-a$, i.e., $c=0$ if $d=a$ or $c=-a$ if $d=-a$.
Because $c^{2}-u'=-u-ac$ by \eqref{1st-8,2},
we have that $u=u'$.

Assume $b\not=0$.
By the equation $4a^{3}-6a^{2}b+4ab^{2}(1-u')+b^{3}(2u'-1)=0$, we have $b$ is even.
Substituting $a=A+\frac{b}{2}$ for some $A\in \Z$ to this equation (i.e., Tschirnhaus's transformation), we have the following equation:
\begin{eqnarray*}
& &4(A+\frac{b}{2})^{3}-6(A+\frac{b}{2})^{2}b+4(A+\frac{b}{2})b^{2}(1-u')+b^{3}(2u'-1) \\
&=&4(A^{3}+3A^{2}\frac{b}{2}+3A\frac{b^{2}}{4}+\frac{b^{3}}{8})-6(A^{2}+Ab+\frac{b^{2}}{4})b+4(Ab^{2}+\frac{b^{3}}{2})(1-u')+b^{3}(2u'-1) \\
&=&4A^{3}+6A^{2}b+3Ab^{2}+\frac{b^{3}}{2}-6A^{2}b-6Ab^{2}-\frac{3b^{3}}{2}+4Ab^{2}+2b^{3}-4Ab^{2}u'-2b^{3}u'+2b^{3}u'-b^{3} \\
&=&4A^{3}+Ab^{2}-4Ab^{2}u' \\
&=&A(4A^{2}+b^{2}-b^{2}u')=0
\end{eqnarray*}
Therefore, there are the two cases: $A=0$ or $A\not=0$.
We first assume $A\not=0$.
Then, by using the equation $4A^{2}+b^{2}-b^{2}u'=0$, we have $u'\ge 1$.
Now, there is the following commutative diagram:
\begin{eqnarray*}
\begin{array}{ccc}
H^{2}(N(u))=\Z X\oplus\Z Y & \stackrel{X}{\longrightarrow} & \Z X^{2}\oplus\Z XY=H^{4}(N(u)) \\
f \downarrow &  & \downarrow f \\
H^{2}(N(u'))=\Z x\oplus\Z y & \stackrel{ax+by}{\longrightarrow} & \Z x^{2}\oplus\Z xy=H^{4}(N(u')) \\
\end{array}
\end{eqnarray*}
Because $X$ and $f$ are isomorphisms, so is $ax+by$ in the diagram. Using the indicated generators as bases,
the determinant of the map $f\circ X:H^{2}(N(u))\to H^{4}(N(u'))$ is equal to the determinant of the map
$(ax+by)\circ f : H^{2}(N(u))\to H^{4}(N(u'))$, which is equal to
\begin{eqnarray}
\label{5th-8,2}
a^{2}-ab+b^{2}u'=\epsilon_{1}=\pm 1.
\end{eqnarray}
Because $a\in \Z$, the discriminant of this equation satisfies
\begin{eqnarray*}
b^{2}-4(b^{2}u'-\epsilon_{1})=b^{2}(1-4u')+4\epsilon_{1}\ge 0 \\
\end{eqnarray*}
Because $u'\ge 1$, we have that
\begin{eqnarray*}
0< b^{2}\le \frac{4\epsilon_{1}}{4u'-1}<1.
\end{eqnarray*}
This gives a contradiction to $b\in \Z$.
Therefore, we have $A=0$, i.e., $a=\frac{b}{2}$.
Because $ad-bc=\epsilon(=\pm 1)$, we also have that $a=\epsilon'=\pm 1$, $b=2\epsilon'$ and $d-2c=\epsilon\epsilon'$.
Hence, by \eqref{5th-8,2},
we have $-1+4u'=\epsilon_{1}$, i.e., $u'=0$ and $\epsilon_{1}=-1$.
By applying a similar method to the one used to derive \eqref{5th-8,2}
for $f^{-1}(x)$, we have
\begin{eqnarray}
\label{6th-8,2}
d^{2}+db+b^{2}u=\epsilon_{2}=\pm 1.
\end{eqnarray}
Substituting \eqref{5th-8,2} and \eqref{6th-8,2} to \eqref{3rd-8,2} and \eqref{4th-8,2},
we have
\begin{eqnarray*}
& & c^{2}=u\epsilon_{1}-u'd^{2}+cd=-u+cd; \\
& & -2ac=\epsilon_{1}+2bdu'-ad-bc=-1-(d+2c)\epsilon'.
\end{eqnarray*}
By using the second equation above, we also have $d=-\epsilon'$; therefore, by $d-2c=\epsilon\epsilon'$, we have $c=\frac{-\epsilon'-\epsilon\epsilon'}{2}=0$ or $-\epsilon'$.
If $c=0$, then $u=0$ by the first equation  above; if $c=-\epsilon'$ then we also have $u=0$ by $d=-\epsilon'$ and the first equation above.
This implies that $u=u'=0$ for the case $b\not=0$.

This establishes the statement.
\end{proof}
Therefore, by Theorem~\ref{Atiyah-Rees} and Lemma \ref{cohom-8dim-3},  we have the following corollary.

\begin{corollary}
\label{N(u)}
The following three statements are equivalent.
\begin{enumerate}
\item Two spaces $N(u)$ and  $N(u')$ are diffeomorphic.
\item Two cohomology rings $H^{*}(N(u))$ and $H^{*}(N(u'))$ are isomorphic.
\item $u=u'\in \Z$.
\end{enumerate}
\end{corollary}

On the other hand, for $M_{\alpha}(u)$ we have the following Proposition.
\begin{proposition}
\label{classification-8dim}
Assume $u(u+1)/12\in \Z$.
The following two statements are equivalent.
\begin{enumerate}
\item Two spaces $M_{\alpha}(u)$ and $M_{\beta}(u')$ are diffeomorphic.
\item $(\alpha,u)=(\beta,u')\in \Z_{2}\times \Z$.
\end{enumerate}
\end{proposition}

In order to prove Proposition~\ref{classification-8dim}, we first compute the $6$-dimensional 
homotopy group of $M_{\alpha}(u)$ in Proposition~\ref{homotopy}.
Now $M_{\alpha}(u)$ can be defined by the following pull-back diagram:
\begin{eqnarray*}
\xymatrix{
& M_{\alpha}(u) \ar[d] \ar[r] & EU(2)\times_{U(2)}\C P^1 \ar[d] \\
& \C P^3 \ar[r]^{\mu_{\alpha,u}}  & BU(2)
}
\end{eqnarray*}

Let $p:S^{7}\to \C P^3$ be the canonical $S^1$-fibration and $P(\xi_{\alpha, u})$ be the pull-back of $M_{\alpha}(u)$ along $p$.
Namely,
we have the following diagram:
\begin{eqnarray}\label{cp1-fibration diagrma}
\xymatrix{
& P(\xi_{\alpha,u}) \ar[d] \ar[r] & M_{\alpha}(u) \ar[d] \ar[r] & EU(2)\times_{U(2)}\C P^1 \ar[d] \\
& S^7 \ar[r]^p & \C P^3 \ar[r]^{\mu_{\alpha,u}}\ar[r]^{\mu_{\alpha,u}} & BU(2) &
}
\end{eqnarray}

Then, we have the following lemma.
\begin{lemma}
\label{isomorphisms}
For $*\ge 3$,
$\pi_{*}(P(\xi_{\alpha,u}))\simeq \pi_{*}(M_{\alpha}(u))$.
\end{lemma}
\begin{proof}
Because $P(\xi_{\alpha,u})$ is the pull-back of $M_{\alpha}(u)$,
the homotopy exact sequences of $P(\xi_{\alpha,u})$ and $M_{\alpha}(u)$ satisfy the following commutative diagram:
\begin{eqnarray*}
\xymatrix{
& \pi_{*+1}(S^7) \ar[r]\ar[d] & \pi_{*}(\C P^1)\ar[r]\ar[d] & \pi_{*}(P(\xi_{\alpha,u}))\ar[r]\ar[d] & \pi_{*}(S^7) \ar[r]\ar[d] & \pi_{*-1}(\C P^1) \ar[d] \\
& \pi_{*+1}(\C P^3) \ar[r] & \pi_{*}(\C P^1)\ar[r] & \pi_{*}(M_{\alpha}(u))\ar[r] & \pi_{*}(\C P^3) \ar[r] & \pi_{*-1}(\C P^1)
}
\end{eqnarray*}
From the homotopy exact sequence of the fibration $S^1\to S^7\to \mathbb CP^3$, we have $\pi_{*}(S^7)\simeq \pi_{*}(\C P^3)$ for $*\ge 3$.
Therefore, by using the $5$ lemma, we have the statement.
\end{proof}

Now we may prove the following proposition.
\begin{proposition}
\label{homotopy}
Assume $u(u+1)/12\in \Z$.
The following two isomorphisms hold.
\begin{enumerate}
\item
$\pi_{6}(P(\xi_{\alpha,u}))\simeq \pi_{6}(M_{\alpha}(u))\simeq \Z_{12}$ if $\alpha\equiv u(u+1)/12\ ({\rm mod}\ 2)$
\item
$\pi_{6}(P(\xi_{\beta,u}))\simeq \pi_{6}(M_{\beta}(u))\simeq \Z_{6}$ if $\beta\not\equiv u(u+1)/12\ ({\rm mod}\ 2)$
\end{enumerate}
\end{proposition}
\begin{proof} 
We first claim the $1$st statement.
If $u(u+1)/12\in \Z$ and $\alpha\equiv u(u+1)/12\  ({\rm mod}\ 2)$, then it follows from \cite{AtRe} that
$\xi_{\alpha,u}$ is induced from the rank $2$ complex vector bundle over $\C P^{4}$.
Namely, there is the following commutative diagram:
\begin{eqnarray}\label{cp1-fibration diagrma}
\xymatrix{
& \xi_{\alpha,u} \ar[d] \ar[r] & \eta_{(\alpha,0,u)} \ar[d] \ar[r] & \widetilde{\mu}_{\alpha,u} \ar[d] \ar[r] & EU(2)\times_{U(2)}\C^2 \ar[d] \\
& S^7 \ar[r]^p & \C P^3 \ar[r]\ar[r] & \C P^{4} \ar[r] & BU(2)
}
\end{eqnarray}
On the other hand, 
we have that  $\pi_{7}(\C P^{4})\simeq \pi_{7}(S^{9})=\{0\}$,
by using the homotopy exact sequence for the fibration $S^{1}\to S^{9}\to \C P^{4}$.
This implies that $\xi_{\alpha,u}$ is the trivial $\C^{2}$-bundle over $S^{7}$.
Therefore, 
\begin{eqnarray*}
P(\xi_{\alpha,u})=S^{7}\times \C P^1
\end{eqnarray*} 
when $u(u+1)/12\in \Z$ and $\alpha\equiv u(u+1)/12\ ({\rm mod}\ 2)$.
Hence, we also have that 
\begin{eqnarray*}
\pi_{6}(M_{\alpha}(u))\simeq \pi_{6}(S^{7}\times \C P^1)\simeq \pi_{6}(\C P^1)\simeq \Z_{12}.
\end{eqnarray*}

Next we claim the $2$nd statement.
Let $\mu_{\alpha,u}:\C P^3\to BU(2)$ be a continuous map which induces the above $\eta_{(\alpha,0,u)}$, and
$\beta$ be the element in $\Z_{2}$ which is not equal to $\alpha$.
Let $x\in \C P^3$ and $s=\mu_{\alpha,u}(x)\in BU(2)$ be base points. 
Take a disk neighborhood around $x\in \C P^3$ and pinch its boundary to a point, i.e., 
the boundary of $D^{6}\subset \C P^3$ pinches to a point, 
then we obtain the surjective map 
\begin{eqnarray*}
\rho:\C P^{3}\to \C P^{3}\vee S^{6},
\end{eqnarray*}
where $\C P^{3}\vee S^{6}$ may be regarded as the wedge sum with respect to the base points $x\in \C P^3$ and $y\in S^{6}$.
Due to theorem of Atiyah-Rees \cite{AtRe}, 
we have
$\eta_{(\beta,0,u)}\not\equiv \eta_{(\alpha,0,u)}$.
This implies that 
the vector bundle $\eta_{(\beta,0,u)}$ is induced from the following continuous map:
\begin{eqnarray}\label{composition}
\xymatrix{
\mu_{\beta,u}:\C P^3 \ar[r]^{\rho} & \C P^{3}\vee S^{6} \ar[r]^{\nu_{\alpha}} & BU(2)
}
\end{eqnarray}
where $\nu_{\alpha}=\mu_{\alpha,u}\vee \kappa$ for the generator $\kappa\in \pi_{6}(BU(2), s)\simeq \Z_{2}$.\footnote{This construction induces the free $\pi_{6}(BU(2))\simeq \pi_{5}(U(2))\simeq \Z_{2}$ action on $\widetilde{KSp}(\C P^{3})\simeq \Z_{2}\oplus\Z$ (see \cite{AtRe}).} 
Hence, we have the following commutative diagram. 
\begin{eqnarray}\label{cp1-fibration diagrma}
\xymatrix{
& P(\xi_{\beta,u}) \ar[d] \ar[r] & M_{\beta}(u) \ar[d] \ar[r] & EU(2)\times_{U(2)}\C P^1 \ar[d] \\
& S^7 \ar[r]^p \ar[rd] & \C P^3 \ar[d]^{\rho}\ar[r]^{\mu_{\beta,u}}\ar[r]^{\mu_{\beta,u}} & BU(2) \\
&    & \C P^{3}\vee S^6 \ar[ru]^{\nu_{\alpha}} &
}
\end{eqnarray}
From the $\C P^1$-fibrations  $\C P^1\to P(\xi_{\beta,u})\to S^7$ and 
 $\C P^1\to EU(2)\times_{U(2)}\C P^1\cong BT^2\to BU(2)$ in the above diagram \eqref{cp1-fibration diagrma},
there is the following commutative diagram.
\begin{eqnarray*}
\xymatrix{
& \pi_{7}(S^7)\simeq \Z \ar[r]\ar[d] & \pi_{6}(\C P^1)\ar[r]\ar[d]^{\simeq} & \pi_{6}(P(\xi_{\beta,u}))\ar[r]\ar[d] & \pi_{6}(S^7)=\{0\}\ar[d] \\
& \pi_{7}(BU(2))\simeq \Z_{12} \ar[r]^{\quad \quad \simeq} & \pi_{6}(\C P^1)\ar[r] & \pi_{6}(BT^2)=\{0\}\ar[r] & \pi_{6}(BU(2))\simeq \Z_{2} 
}
\end{eqnarray*}
This diagram shows that the following exact sequence:
\begin{eqnarray}
\label{final-ex}
\Z\simeq \pi_{7}(S^{7})\to \pi_{7}(BU(2))(\simeq\Z_{12})\to \pi_{6}(P(\xi_{\beta,u}))\to \{0\}.
\end{eqnarray}
In this diagram, the left homomorphism is induced from $\widetilde{\mu}:=\mu_{\beta,u}\circ p \colon S^7\to BU(2)$, 
say $\widetilde{\mu}_{\#}:\Z\to\Z_{12}$.
We claim $\widetilde{\mu}_{\#}(1)=[6]_{12}\in \Z_{12}$.
Because the diagram \eqref{cp1-fibration diagrma} is commutative, we may regard that 
$\widetilde{\mu}:=\mu_{\beta,u}\circ p \colon S^7\to BU(2)$ can be defined by passing through the map $\nu_{\alpha}:\C P^{3}\vee S^6\to BU(2)$,  i.e.,
$\widetilde{\mu}=\nu_{\alpha}\circ\rho\circ p$.
Because $\nu_{\alpha}=\mu_{\alpha,u}\vee \kappa$,
we also have 
\begin{eqnarray*}
\widetilde{\mu}=(\mu_{\alpha,u}\vee \kappa)\circ\rho\circ p=(\mu_{\alpha,u}\circ\rho\circ p)\vee (\kappa\circ\rho\circ p).
\end{eqnarray*}
By the argument when we proved the $1$st statement, we see that 
$\mu_{\alpha,u}\circ\rho\circ p$ induces the trivial bundle over $S^{7}$, i.e., $\mu_{\alpha,u}\circ\rho\circ p$ is homotopic to the trivial map.
This also implies that 
there is the following decomposition up to homotopy:
\begin{eqnarray*}
\widetilde{\mu}:S^{7}\stackrel{p}{\longrightarrow} \C P^{3}\stackrel{\rho}{\longrightarrow} \C P^{3}\vee S^{6}\stackrel{\pi}{\longrightarrow} S^{6}\stackrel{\kappa}{\longrightarrow} BU(2),
\end{eqnarray*}
where $\pi$ is the collapsing map of $\C P^3$ to a point.
Therefore, we have the following decomposition for the induced map
\begin{eqnarray*}
\widetilde{\mu}_{\#}:\pi_{7}(S^{7})\stackrel{\Psi_{\#}}{\longrightarrow} \pi_{7}(S^{6})\simeq \Z_{2}\stackrel{\kappa_{\#}}{\longrightarrow} \pi_{7}(BU(2))\simeq \Z_{12},
\end{eqnarray*}
where the $1$st map is induced from the surjective map $\Psi=\pi\circ \rho\circ p$.
Because $\Psi$ is non trivial map, $\Psi_{\#}(1)=[1]_{2}$ (the generator of $\pi_{7}(S^{6})\simeq \Z_{2}$).
Moreover, because $\kappa\in \pi_{6}(BU(2))\simeq \Z_{2}$ is the generator, i.e., non-trivial map,
we have $\kappa_{\#}([1]_{2})=[6]_{12}\in \Z_{12}$.
This shows that $\widetilde{\mu}_{\#}(1)=[6]_{12}$; therefore, $\widetilde{\mu}_{\#}(\pi_{7}(S^{7}))=\{[0]_{12},\ [6]_{12}\}\subset \Z_{12}$.

Consequently, by the exact sequence \eqref{final-ex}, we have that 
\begin{eqnarray*}
\pi_{6}(P(\xi_{\beta,u}))\simeq \pi_{7}(BU(2))/\widetilde{\mu}_{\#}(\pi_{7}(S^{7}))\simeq \Z_{12}/\{[0]_{12}, [6]_{12}\}\simeq \Z_{6}.
\end{eqnarray*}
By Lemma \ref{isomorphisms}, we have the statement.
\end{proof}

\begin{remark}
For example, the relation $u(u+1)/12\in \Z$ is true for the case when $u=0$ and $u=3$.
In these cases, by using Proposition \ref{homotopy}, we have
\begin{eqnarray*}
\pi_{6}(M_{\alpha}(0))\simeq
\left\{
\begin{array}{cc}
\Z_{12} & {\rm for}\ \alpha\equiv 0 \\
\Z_{6} & {\rm for}\ \alpha\equiv 1
\end{array}
\right.
\end{eqnarray*}
and
\begin{eqnarray*}
\pi_{6}(M_{\alpha}(3))\simeq
\left\{
\begin{array}{cc}
\Z_{6} & {\rm for}\ \alpha\equiv 0 \\
\Z_{12} & {\rm for}\ \alpha\equiv 1
\end{array}
\right.
\end{eqnarray*}

On the other hand, the case when $u=1$ does not satisfy the relation $u(u+1)/12\in \Z$.
It follows from the cohomology ring of the flag manifold of type $C$ (see e.g. \cite{Bo} or \cite{FIM}) 
that the flag manifold $Sp(2)/T^{2}$ is one of this case, i.e., 
$M_{0}(1)$ or $M_{1}(1)$.
However, by using the homotopy exact sequence for the fibration $T^{2}\to Sp(2)\to Sp(2)/T^{2}$ and the computation in \cite{MiTo}, 
we have that
\begin{eqnarray*}
\pi_{6}(Sp(2)/T^{2})\simeq \pi_{6}(Sp(2))=0.
\end{eqnarray*} 
Therefore, Proposition \ref{homotopy} is not true for the case when $u(u+1)/12\not\in \Z$.
\end{remark}

Let us prove Proposition \ref{classification-8dim}

\begin{proof}[Proof of Proposition \ref{classification-8dim}]
By using Theorem \ref{Atiyah-Rees}, $(2)\Rightarrow (1)$ is trivial.
We prove $(1)\Rightarrow (2)$.
Assume $M_{\alpha}(u)\cong M_{\beta}(u')$.
If $u\not=u'$, then $H^{*}(M_{\alpha}(u))\not\simeq H^{*}(M_{\beta}(u'))$ by Lemma \ref{cohom-8dim-2}.
Therefore, we have $u=u'$.
By Proposition \ref{homotopy}, $M_{0}(u)\not\cong M_{1}(u)$.
This implies that the statement $(1)\Rightarrow (2)$ in Proposition \ref{classification-8dim}.
This establishes Proposition \ref{classification-8dim}.
\end{proof}

Consequently, by Lemma \ref{cohom-8dim-2} and Proposition \ref{classification-8dim}, we have the following corollary:
\begin{corollary}
The set of $8$-dimensional $\C P$-towers does not satisfy the cohomological rigidity.
\end{corollary}
Note that if we restrict the class of $8$-dimensional $\C P$-towers to the $8$-dimensional generalized Bott manifolds with height $2$,
then cohomological rigidity holds by \cite{CMS10}.

Using Corollary \ref{N(u)} and Proposition \ref{classification-8dim}, we also have Theorem \ref{main-thm2}.

\section*{Acknowledgments}
The authors would like to thank Mikiya Masuda for his helpful comments about torus actions on $\C P$-towers.
The first author would like to give heartful thanks to Prof.\ Nigel Ray whose comments and helps to stay in University of Manchester
were innumerably valuable. 
He would also like to thank Takahiko Yoshida and Shizuo Kaji whose useful comments were an enormous help to him.

\end{document}